\documentclass[review,hidelinks,onefignum,onetabnum]{siamart251216}
\newtheorem{assumption}{Assumption}

\usepackage{adjustbox}
\usepackage{booktabs}

\usepackage{lipsum}
\usepackage{amsfonts}
\usepackage{graphicx}
\usepackage{epstopdf}
\usepackage{algorithmic}
\ifpdf
  \DeclareGraphicsExtensions{.eps,.pdf,.png,.jpg}
\else
  \DeclareGraphicsExtensions{.eps}
\fi

\newsiamremark{remark}{Remark}
\newsiamremark{hypothesis}{Hypothesis}
\crefname{hypothesis}{Hypothesis}{Hypotheses}
\newsiamthm{claim}{Claim}
\newsiamremark{fact}{Fact}
\crefname{fact}{Fact}{Facts}

\headers{DSAN for Elliptic Interface Problems}{J. Luo, H. Kang, X. Yue and S. Xu}

\title{A Domain-Specific Adaptation Network for Elliptic Interface Problems\thanks{The work is supported by the National Natural Science Foundation of China
(Nos. 12671501 and 12571409), the National Key R\&D Program of China
(No. 2024YFA1016300), the Basic Research Program of Jiangsu
(No. BK20250803), the Key University Science Research Project of Jiangsu
Province (No. 25KJB110014), and the Collaborative Team Project of the Open
Research Program of the National Key Laboratory of Mathematical Sciences
(No. SLMS-2025-KFKT-TD-04).}}

\author{Jiaqi Luo\thanks{School of Mathematical Sciences, Soochow University, Suzhou, Jiangsu, China
  (\email{jqluo@suda.edu.cn},\email{khm@suda.edu.cn},\email{xyyue@suda.edu.cn}).}
\and Hongmei Kang\footnotemark[2]
\and Xingye Yue\footnotemark[2] 
\and Shixin Xu\thanks{Zu Chongzhi Center, Duke Kunshan University, Kunshan, Jiangsu, China
  (Corresponding auther, \email{shixin.xu@dukekunshan.edu.cn})}
}

\usepackage{amsopn}

\ifpdf
\hypersetup{
  pdftitle={A Domain-Specific Adaptation Network for Elliptic Interface Problems},
  pdfauthor={J. Luo, H. Kang, X. Yue and S. Xu}
}
\fi

\allowdisplaybreaks \allowdisplaybreaks[4]
\begin{document}

\maketitle

\begin{abstract}
Deep learning provides a flexible mesh-free framework for solving elliptic interface problems. However, domain-specific neural approaches often rely on separate networks for individual subdomains, causing the model complexity to grow rapidly as the number of subdomains increases. To address this limitation, we propose a \textbf{Domain-Specific Adaptation Network (DSAN)}, which combines a shared backbone with lightweight subdomain-specific adapters to provide local representation flexibility without replicating a complete network for each subdomain. We further establish a stability and error analysis showing that the continuous physics-informed residual controls the solution error, and decompose the resulting error into approximation, optimization, and sampling components. Numerical experiments on two- and three-dimensional problems with multiple regular and irregular interfaces demonstrate that DSAN achieves accuracy comparable to or better than fully decomposed multi-network methods while requiring substantially fewer trainable parameters. Moreover, the parameter savings become increasingly pronounced as the number of subdomains grows, demonstrating the scalability of the proposed domain-specific adaptation framework.
\end{abstract}



\section{Introduction}
Interface problems arise in a wide range of scientific and engineering applications, including material science \cite{medlin2009interfaces}, fluid dynamics \cite{jin2021nsfnets}, and heat transfer \cite{cai2021physics}. They are typically characterized by non-smooth behavior across interfaces that separate the computational domain into multiple subdomains, where the governing equations, coefficients, or interface conditions may differ. Compared with standard partial differential equations (PDEs), interface problems are more difficult to solve numerically because the solution must satisfy not only the governing equations within individual subdomains but also the prescribed jump and flux conditions across interfaces.

A variety of classical numerical methods have been developed for interface problems, including the extended finite element method \cite{xiao2020high}, level-set methods \cite{sussman1994level}, immersed boundary methods \cite{peskin2002immersed}, and immersed interface methods \cite{leveque1994immersed,li2006immersed}. These approaches have established mature numerical frameworks and exhibit good performance for many problems. However, their practical implementation often depends on mesh construction or geometry-aware discretization. When the interfaces become geometrically complicated, mesh generation and interface treatment can become increasingly difficult.

Recently, deep learning-based PDE solvers have emerged as mesh-free alternatives to conventional numerical methods \cite{raissi2019physics,weinan2018deep,han2018solving,zang2020weak,weng2026deep}. 
These methods approximate the solution directly from the governing equations and associated constraints, constructing a discrete objective function using randomly sampled points~\cite{wu2023comprehensive,luo2025imbalanced,gao2023failure}. By avoiding explicit mesh generation, they are particularly attractive for problems with complex geometries.

A common strategy is to retain a \textit{global representation}, in which a unified neural model is used over the entire computational domain while additional information is introduced to distinguish different solution branches. Hu et al.~\cite{hu2022discontinuity} proposed the discontinuity capturing shallow neural network (DCSNN), where an additional coordinate is used to label different subdomains, thereby extending a piecewise continuous function to a continuous representation in a higher-dimensional space. More recently, a categorical embedding strategy was introduced for anisotropic interface problems with many subdomains~\cite{hu2026categorical}, where one-hot-encoded subdomain information is mapped to a lower-dimensional embedding to improve the scalability of feature augmentation.
Related ideas have also been combined with extreme learning machine architectures~\cite{song2026discontinuity,lin2025discontinuous}, for which training can be reduced to a least-squares problem. Other approaches incorporate level-set-based feature augmentation into the network input to encode interface geometry and non-smooth solution behavior ~\cite{zhang2025immersed,tseng2023cusp,lai2022shallow}, while \cite{fan2025novel,li2025two} construct global representations by decomposing the interface solution into continuous and discontinuous components. These methods preserve a compact single-network formulation by encoding subdomain dependence through augmented inputs or solution decomposition.

For challenging interface problems, however, a unified global representation may provide limited flexibility for capturing strongly
different behaviors across subdomains. This has motivated a second line of work based on \textit{domain-specific representations}, where different subdomains are assigned distinct representation mechanisms and are coupled through interface conditions.
The Extended Physics-Informed Neural Network (XPINN)~\cite{jagtap2020extended} is an early method that assigns an individual neural network to each subdomain and couples the local solutions through interface conditions. 
Subsequent studies have explored alternative network constructions to improve approximation accuracy or training efficiency, including transferable networks~\cite{lu2025multiple}, local randomized neural networks~\cite{li2025local}, and the random feature  method~\cite{chi2024random}. 
In \cite{he2022mesh}, a classification mechanism is introduced to determine the subdomain associated with each sampling point and to select the corresponding local network, providing greater flexibility for approximating solutions with strong singular or discontinuous behavior near interfaces. 
The interfaced neural network (INN)~\cite{wu2022inn} also adopts a multi-network formulation and introduces an extended multi-gradient descent strategy to adaptively balance different components of the training objective. Interface-PINN (I-PINN)~\cite{sarma2024interface} introduces subdomain-specific activations while sharing the remaining network parameters across regions, thereby providing domain-dependent representations without fully duplicating the model.
IG-PINN~\cite{zheng2025ig} combines a global fully connected network that captures the smooth solution component with interface-gated local networks that provide subdomain-specific corrections.
By introducing explicit regional specialization, these methods generally offer greater representation flexibility and often achieve stronger approximation performance on difficult interface problems.

Domain-specific formulations are attractive because they provide dedicated representation flexibility for individual subdomains and can effectively capture discontinuous or strongly heterogeneous solution behavior. However, this flexibility is often achieved by assigning a separate neural network to each subdomain. As the number of interfaces increases, the resulting growth in model size can become substantial, limiting the scalability of such approaches to problems with many subdomains.

To address this issue, we propose a \textit{Domain-Specific Adaptation Network} (DSAN) for elliptic interface problems. DSAN employs a shared backbone across the computational domain and introduces lightweight adapters for individual subdomains. In this way, the main network parameters are shared, while subdomain-dependent degrees of freedom are provided through the adapters. This design retains the flexibility of domain-specific representations without replicating a complete network for each subdomain, thereby reducing the growth in trainable parameters as the number of subdomains increases. We further establish a stability and error analysis for the physics-informed formulation, characterize the approximation, optimization, and sampling contributions to the solution error, and derive the parameter complexity of DSAN. 
Numerical experiments on two- and three-dimensional problems with multiple regular and irregular interfaces demonstrate that DSAN achieves accuracy comparable to or better than fully decomposed multi-network methods while using substantially fewer trainable parameters.

The remainder of this paper is organized as follows. Section~\ref{s.method} introduces the elliptic interface problem and presents the DSAN architecture and physics-informed training formulation. Section~\ref{s.erroranalysis} develops the stability and error analysis and studies the parameter complexity of the proposed method. Section~\ref{s.exps} presents the numerical experiments and examines the scalability of DSAN with respect to the number of interfaces. Section~\ref{s:conclu} concludes the paper.





\section{Methodology}
\label{s.method}
\subsection{Elliptic Interface Problems}
Let $\{\Gamma_k\}_{k=1}^{K}$ denote a collection of mutually non-intersecting interfaces contained in $\Omega$. These interfaces partition $\Omega$ into $K+1$ non-overlapping connected subdomains, denoted by $\{\Omega_j\}_{j=0}^{K}$, such that
\begin{equation}
\Omega\setminus\bigcup_{k=1}^{K}\Gamma_k=\bigcup_{j=0}^{K}\Omega_j,~\Omega_i\cap\Omega_j=\emptyset~\text{for }i\neq j.
\end{equation}
We designate $\Omega_0$ as the background subdomain adjacent to the outer boundary $\partial\Omega$.

For each interface $\Gamma_k$, let $i_k^{+}$ and $i_k^{-}$ denote the indices of the two subdomains adjacent to $\Gamma_k$. We denote these subdomains by $\Omega_{i_k^{+}}$ and $\Omega_{i_k^{-}}$, respectively, and let $\mathbf n_k$ be the unit normal pointing from $\Omega_{i_k^{+}}$ to $\Omega_{i_k^{-}}$.

The governing equations are given by
\begin{subequations}
\label{e.interface}
\begin{align}
-\nabla \cdot \bigl(\beta(\mathbf{x}) \nabla u(\mathbf{x})\bigr)
&= f(\mathbf{x}),&& \mathbf{x} \in \Omega \setminus \bigcup_{k=1}^{K} \Gamma_k,
\label{e.interface.a} \\[4pt]
[u]_{\Gamma_k}
&= v_k(\mathbf{x}),&& \mathbf{x} \in \Gamma_k,~k=1,\ldots,K,
\label{e.interface.b} \\[4pt]
[\beta\nabla u\cdot\mathbf n_k]_{\Gamma_k}
&= w_k(\mathbf{x}),
&& \mathbf{x} \in \Gamma_k,~k=1,\ldots,K,
\label{e.interface.c} \\[4pt]
u(\mathbf{x})
&= g(\mathbf{x}),&& \mathbf{x} \in \partial\Omega.
\label{e.interface.d}
\end{align}
\end{subequations}
\noindent
Here, $\mathbf{x}=(x_1,\ldots,x_d)\in\mathbb{R}^d$ denotes the spatial variable, while $\nabla$ and $\nabla\cdot$ denote the gradient and divergence operators, respectively.  The coefficient $\beta(\mathbf{x})\in\mathbb{R}^{d\times d}$ is assumed to be bounded and uniformly positive definite. 
Specifically, there exist constants $0<\beta_{\min}\le \beta_{\max}<\infty$ such that
\begin{equation}\label{ass:uniformellipticity}
\beta_{\min}|\boldsymbol{\xi}|^2\le\boldsymbol{\xi}^{T}\beta(\mathbf{x})\boldsymbol{\xi}\le
\beta_{\max}|\boldsymbol{\xi}|^2,\qquad\forall\,\boldsymbol{\xi}\in\mathbb{R}^d,
\end{equation}
\noindent
for almost every $\mathbf{x}\in\Omega$. 
The coefficient may be discontinuous across the interfaces but is assumed to be sufficiently smooth within each subdomain. The functions $f$ and $g$ denote the source term and the prescribed Dirichlet boundary data, respectively.

For any quantity $q$ admitting one-sided traces on $\Gamma_k$, let $q^{+}$ and $q^{-}$ denote the traces taken from $\Omega_{i_k^{+}}$ and $\Omega_{i_k^{-}}$, respectively. The oriented jump across $\Gamma_k$ is
defined by $[q]_{\Gamma_k}=q^{+}-q^{-}$.
Accordingly, $v_k$ and $w_k$ prescribe the jumps in the solution and the
normal flux across $\Gamma_k$, respectively.

When $f$ and $\beta$ are sufficiently smooth within each subdomain, the solution $u$ is generally piecewise smooth, while its value and derivatives may be discontinuous across the interfaces. Consequently, a neural-network approximation of such interface problems requires an appropriate mechanism to represent and distinguish the solution behavior in different subdomains.

\subsection{Domain-Specific Representation}

A direct strategy for representing the piecewise solution of
Eq.~\eqref{e.interface} is to construct a domain-specific approximation.
For each subdomain $\Omega_j$, $j=0,\ldots,K$, the local approximation is
represented as
\begin{equation}
u_{\theta,j}(\mathbf{x}) = \mathcal{N}_j(\mathbf{x};\Theta_j), ~\mathbf{x}\in\Omega_j,
\end{equation}
where $\mathcal{N}_j$ denotes the neural network associated with the
$j$-th subdomain and $\Theta_j$ denotes its parameters. The complete
approximation is consequently defined in a piecewise manner as
\begin{equation}
u_\theta(\mathbf{x})
=
u_{\theta,j}(\mathbf{x}),
\qquad
\mathbf{x}\in\Omega_j
\end{equation}

The local representations are coupled through the interface conditions in Eq.~\eqref{e.interface}. At an interface $\Gamma_k$, the two one-sided predictions associated with the neighboring subdomains must satisfy
\begin{equation}
\begin{aligned}
u_{\theta,i_k^{+}}(\mathbf x)-u_{\theta,i_k^{-}}(\mathbf x)
&=v_k(\mathbf x), && \mathbf x\in\Gamma_k, \\
\left(\beta_{i_k^{+}}\nabla u_{\theta,i_k^{+}}
-\beta_{i_k^{-}}\nabla u_{\theta,i_k^{-}}\right)\cdot\mathbf n_k
&=w_k(\mathbf x), && \mathbf x\in\Gamma_k.
\end{aligned}
\end{equation}
Therefore, domain-specific representation provides separate approximation flexibility within  individual regions while maintaining consistency
between neighboring solutions through the interface conditions.

A conventional implementation of this strategy assigns a fully independent neural network to each subdomain, with no parameters shared between networks. For $K+1$ subdomains, this requires $K+1$ separate sets of network parameters. Consequently, for a fixed local network architecture, the total parameter count grows linearly with the number of subdomains.

\subsection{Domain-Specific Adaptation Network for Elliptic Interface Problems}
\label{ss.network}
\begin{figure}[!ht]
    \centering
    \includegraphics[width=0.95\linewidth]{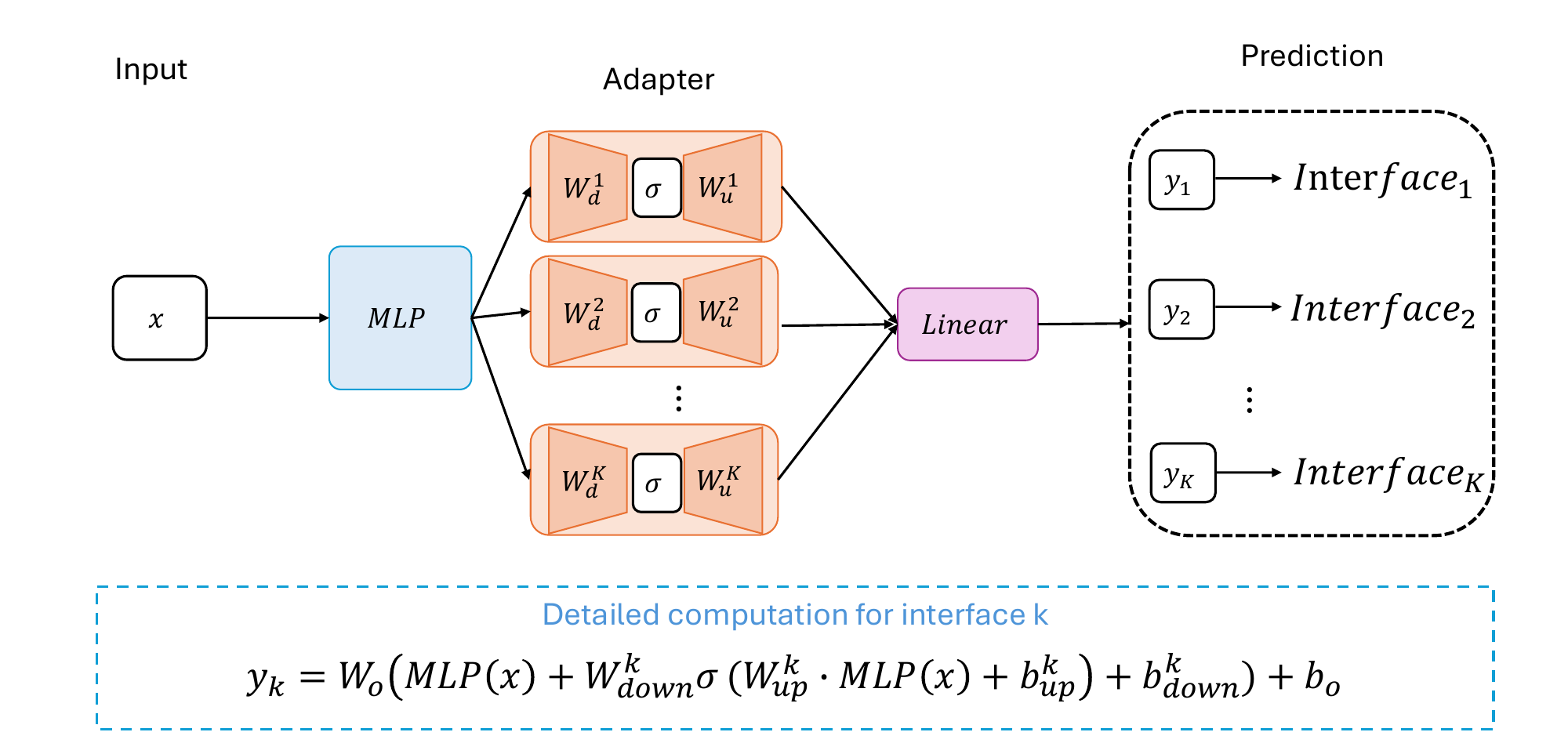}
    \caption{Architecture of the proposed Domain-Specific Adaptation Network.}
    \label{f.network}
\end{figure}

Although the governing differential operator is shared across all subdomains, the corresponding local solutions may differ substantially due to variations in geometric location, interface conditions, and surrounding geometric configuration. From a learning perspective, the solution in each subdomain can therefore be viewed as belonging to a distinct domain-specific distribution, while retaining common structure imposed by the same underlying physical law. This naturally motivates a domain-adaptation perspective, in which shared components capture common physical features across subdomains, while domain-specific components account for variations in local solution behavior.

Based on this motivation, we propose a \textit{Domain-Specific Adaptation Network} (DSAN). The model uses a shared backbone to construct a common representation and introduces a learnable adaptation module for each subdomain to modify this representation according to its local requirements. This architecture provides domain-specific approximation flexibility while reducing parameter duplication compared with fully independent networks.

\paragraph{Network Architecture}

As illustrated in Fig.~\ref{f.network}, DSAN consists of a shared backbone, subdomain-specific adapters, and a shared output layer. 
For an input point $\mathbf{x}\in\Omega_k$, the shared backbone first maps the spatial coordinate to a latent representation,
\begin{equation}
\mathbf{h}=B(\mathbf{x};\Theta_{\mathrm b}),~\mathbf{h}\in\mathbb{R}^{D},
\end{equation}
where $B(\cdot;\Theta_{\mathrm b})$ denotes a multi-layer perceptron (MLP)  with parameters $\Theta_{\mathrm b}$ shared across all subdomains.

Subdomain-specific adaptation is then introduced through a lightweight  bottleneck adapter,
\begin{equation}
\mathbf{h}_j=\mathbf{h}+A_j(\mathbf{h};\Phi_j),~j=0,\ldots,K,
\label{e.adapter}
\end{equation}
where $A_j(\mathbf{h};\Phi_j)=W_{\mathrm{up}}^j\sigma\left(W_{\mathrm{down}}^j\mathbf{h}+\mathbf{b}_{\mathrm{down}}^j\right)+\mathbf{b}_{\mathrm{up}}^j$.
Here, $W_{\mathrm{down}}^j\in\mathbb{R}^{r\times D},
W_{\mathrm{up}}^j\in\mathbb{R}^{D\times r},
\mathbf{b}_{\mathrm{down}}^j\in\mathbb{R}^{r},
\mathbf{b}_{\mathrm{up}}^j\in\mathbb{R}^{D}$,
with $r\ll D$, and $\Phi_j$ denotes the collection of trainable adapter parameters associated with subdomain $\Omega_j$. 

The residual form in \eqref{e.adapter} retains the shared latent representation while introducing a low-dimensional, learnable correction that captures subdomain-specific variations in the local solution. 
Unlike conventional adapter-based model adaptation~\cite{houlsby2019parameter}, where adapters are introduced to adapt a pretrained model to different tasks or domains, here the adapters are embedded into a single PDE solution representation and are indexed by physical subdomains. Consequently, different one-sided solution branches can be represented through lightweight regional components while the main network parameters remain shared. 

Finally, a shared output layer maps the adapted representation to the local solution. For any $\mathbf{x}\in\Omega_j$,
\begin{align}
u_{\theta,j}(\mathbf{x})&=\mathcal{O}(\mathbf{h}_j;\Theta_{\mathrm o})=W_{\mathrm o}\mathbf{h}_j+b_{\mathrm o}
= u_{\theta}^{\mathrm{sh}}(\mathbf{x}) 
+ u_{\theta,j}^{\mathrm{ad}}(\mathbf{x}).
\label{e.dsan_decomposition}
\end{align}
Here, \(u_{\theta}^{\mathrm{sh}}(\mathbf{x})=W_{\mathrm o}B(\mathbf{x};\Theta_{\mathrm b})+b_{\mathrm o}\)
denotes the shared component common to all subdomains, while
\(u_{\theta,j}^{\mathrm{ad}}(\mathbf{x})=W_{\mathrm o}A_j\left(B(\mathbf{x};\Theta_{\mathrm b});\Phi_j\right)\)
represents the subdomain-specific adaptation associated with $\Omega_j$.

\paragraph{Interface Evaluation}

At an interface, the solution admits two one-sided representations associated with the neighboring subdomains. Consider an interface $\Gamma_k$ separating
$\Omega_{i_k^{+}}$ and $\Omega_{i_k^{-}}$. For an interface point
$\mathbf{x}_{\Gamma}\in\Gamma_k$, the shared backbone first produces
\begin{equation}
\mathbf h_{\Gamma}=B(\mathbf{x}_{\Gamma};\Theta_{\mathrm b}).
\end{equation}
The same shared representation is then processed by the adapters associated
with the two neighboring subdomains,
\begin{equation}
\mathbf h_k^{+}=\mathbf h_{\Gamma}+A_{i_k^{+}}(\mathbf h_{\Gamma};\Phi_{i_k^{+}}),~\mathbf h_k^{-}=\mathbf h_{\Gamma}+A_{i_k^{-}}(\mathbf h_{\Gamma};\Phi_{i_k^{-}}).
\end{equation}
The corresponding one-sided predictions are
\begin{equation}
u_{\theta,k}^{+}(\mathbf{x}_{\Gamma})=W_{\mathrm o}\mathbf h_k^{+}+b_{\mathrm o},~
u_{\theta,k}^{-}(\mathbf{x}_{\Gamma})=W_{\mathrm o}\mathbf h_k^{-}+b_{\mathrm o}.
\end{equation}

Although the two predictions are evaluated at the same spatial coordinate, they may take different values because different subdomain-specific adapters are applied. The predicted solution jump is therefore
\begin{equation}
\label{e.jumpnet}
[u_\theta]_{\Gamma_k}=u_{\theta,k}^{+}-u_{\theta,k}^{-}.
\end{equation}
Similarly, the predicted normal-flux jump is
\begin{equation}
\label{e.fluxnet}
[\beta\nabla u_\theta\cdot\mathbf n_k]_{\Gamma_k}=\left(\beta_{i_k^{+}}\nabla u_{\theta,k}^{+}-\beta_{i_k^{-}}\nabla u_{\theta,k}^{-}\right)\cdot\mathbf n_k .
\end{equation}

\paragraph{Physics-Informed Training Objective}

The DSAN is trained by minimizing a physics-informed loss that enforces the governing equation in each subdomain together with the boundary and interface conditions,
\begin{equation}
\mathcal{L}=\lambda_{\mathrm{PDE}}\mathcal{L}_{\mathrm{PDE}}+
\lambda_{\mathrm{BC}}\mathcal{L}_{\mathrm{BC}}+
\lambda_{\mathrm{jump}}\mathcal{L}_{\mathrm{jump}}+
\lambda_{\mathrm{flux}}\mathcal{L}_{\mathrm{flux}}.
\label{e.total_loss}
\end{equation}

Let $\{\mathbf{x}_{j,i}\}_{i=1}^{N_j}$ denote the collocation points sampled from subdomain $\Omega_j$, for $j=0,\ldots,K$. The PDE residual is evaluated separately in each subdomain using the corresponding domain-specific
prediction $u_j$,
\begin{equation}
\mathcal{L}_{\mathrm{PDE}}=\frac{1}{K+1}\sum_{j=0}^{K}\frac{1}{N_j}\sum_{i=1}^{N_j}
\left|-\nabla\cdot\left(\beta_j(\mathbf{x}_{j,i})\nabla u_j(\mathbf{x}_{j,i})\right)-f_j(\mathbf{x}_{j,i})\right|^2,
\label{e.loss_pde}
\end{equation}
where $\beta_j$ and $f_j$ denote the restrictions of $\beta$ and $f$ to
$\Omega_j$, respectively.

For boundary collocation points $\{\mathbf{x}_{b,i}\}_{i=1}^{N_b}\subset\partial\Omega$, the Dirichlet boundary condition is enforced through
\begin{equation}
\mathcal{L}_{\mathrm{BC}}=\frac{1}{N_b}\sum_{i=1}^{N_b}
\left|u_0(\mathbf{x}_{b,i})-g(\mathbf{x}_{b,i})\right|^2,
\label{e.loss_bc}
\end{equation}
where $u_0$ is used because the outer boundary $\partial\Omega$ is adjacent to the background subdomain $\Omega_0$.

For each interface $\Gamma_k$, let $\{\mathbf{x}_{k,i}^{\Gamma}\}_{i=1}^{N_{\Gamma_k}}$ denote the corresponding interface collocation points. The solution-jump condition is enforced by
\begin{equation}
\mathcal{L}_{\mathrm{jump}}=\frac{1}{K}\sum_{k=1}^{K}\frac{1}{N_{\Gamma_k}}\sum_{i=1}^{N_{\Gamma_k}}
\left|
[u_\theta]_{\Gamma_k}(\mathbf{x}_{k,i}^{\Gamma})-v_k(\mathbf{x}_{k,i}^{\Gamma})
\right|^2,
\label{e.loss_jump}
\end{equation}
where the jump is evaluated according to the convention introduced in \eqref{e.jumpnet}.
Similarly, the normal-flux jump condition is enforced through
\begin{equation}
\mathcal{L}_{\mathrm{flux}}=\frac{1}{K}\sum_{k=1}^{K}\frac{1}{N_{\Gamma_k}}\sum_{i=1}^{N_{\Gamma_k}}
\left|
[\beta\nabla u_\theta\cdot\mathbf n_k]_{\Gamma_k}(\mathbf{x}_{k,i}^{\Gamma})-w_k(\mathbf{x}_{k,i}^{\Gamma})
\right|^2,
\label{e.loss_flux}
\end{equation}
where the flux jump is evaluated according to \eqref{e.fluxnet}.

The coefficients $\lambda_{\mathrm{PDE}}$, $\lambda_{\mathrm{BC}}$, $\lambda_{\mathrm{jump}}$, and $\lambda_{\mathrm{flux}}$ in \eqref{e.total_loss} control the relative contributions of the corresponding loss components. All shared parameters, $\Theta_{\mathrm b}$ and $\Theta_{\mathrm o}$, together with the subdomain-specific adapter parameters $\{\Phi_j\}_{j=0}^{K}$, are optimized jointly under the objective \eqref{e.total_loss}.

\section{Stability, Error and Complexity Analysis}
\label{s.erroranalysis}

This section develops the theoretical analysis of DSAN in terms of both solution accuracy and parameter efficiency. We first establish a deterministic residual-to-solution stability estimate, then relate the empirical loss to its continuous population counterpart, and finally combine the stability, sampling, optimization, and approximation estimates to obtain a total error bound for the trained DSAN. The approximation term is further interpreted through the shared-plus-adaptation architecture. We subsequently analyze the parameter complexity of DSAN and its scaling with the number of subdomains.

Throughout this section, we assume that the elliptic interface problem satisfies the following regularity conditions.

\begin{assumption}
\label{ass:regularity}
The computational domain $\Omega\subset\mathbb R^d$ is bounded with $C^{1,1}$ boundary. The interfaces $\{\Gamma_k\}_{k=1}^K$ are mutually disjoint $C^{1,1}$ hypersurfaces contained in $\Omega$. Each interface $\Gamma_k$ separates two neighboring subdomains $\Omega_{i_k^{+}}$ and $\Omega_{i_k^{-}}$, and the prescribed unit normal $\mathbf n_k$ points from $\Omega_{i_k^{+}}$ to $\Omega_{i_k^{-}}$.

For each $j=0,\ldots,K$, the coefficient $\beta_j=\beta|_{\Omega_j}$ is symmetric and satisfies
$\beta_j\in W^{1,\infty}(\Omega_j;\mathbb R^{d\times d})$,
together with the uniform ellipticity condition \eqref{ass:uniformellipticity}.
The problem data satisfy $f_j\in L^2(\Omega_j),~g\in H^{3/2}(\partial\Omega),~v_k\in H^{3/2}(\Gamma_k),~w_k\in H^{1/2}(\Gamma_k)$,
for $j=0,\ldots,K$ and $k=1,\ldots,K$.
\end{assumption}



These assumptions provide a convenient sufficient setting for the analysis. More general geometries and coefficients can also be considered provided that the corresponding elliptic transmission problem satisfies the piecewise $H^2$ regularity estimate used below.
To measure the approximation error, let $q=\{q_j\}_{j=0}^K$ be a piecewise-defined function with
$q_j=q|_{\Omega_j}$. We use $\|q\|_{H^2_{\mathrm{pw}}(\Omega)}^2=\sum_{j=0}^K\|q_j\|_{H^2(\Omega_j)}^2,~
\|q\|_{L^2(\Omega)}^2=\sum_{j=0}^K\|q_j\|_{L^2(\Omega_j)}^2$.
For clarity of presentation, the proofs of all theoretical results in this section are collected in the Appendix \ref{a.proof}.

\subsection{Residual Stability Analysis}
\label{ss.residual}
We first establish a deterministic stability estimate showing that small continuous PDE, boundary, and interface residuals imply a small solution error. The argument relies on the standard piecewise $H^2$ regularity of elliptic transmission problems, which we state below.

\begin{lemma}[Piecewise elliptic regularity]
\label{lem:elliptic_regularity}
Under Assumption~\ref{ass:regularity}, the elliptic transmission problem
\eqref{e.interface} admits a unique piecewise $H^2$ solution.
Moreover, for any
$z=\{z_j\}_{j=0}^K$ with $z_j\in L^2(\Omega_j)$, the homogeneous
transmission problem
\begin{subequations}
\label{e.z}
\begin{align}
-\nabla\cdot(\beta_j\nabla\phi_j)
&=z_j,&& \mathbf x\in\Omega_j,~ j=0,\ldots,K,\\
\phi_{i_k^+}-\phi_{i_k^-}
&=0,&& \mathbf x\in\Gamma_k,~k=1,\ldots,K,\\
\left(
\beta_{i_k^+}\nabla\phi_{i_k^+}-\beta_{i_k^-}\nabla\phi_{i_k^-}\right)\cdot\mathbf n_k
&=0,&& \mathbf x\in\Gamma_k,~k=1,\ldots,K,\\
\phi_0&=0,&& \mathbf x\in\partial\Omega,
\end{align}
\end{subequations}
has a unique piecewise $H^2$ solution $\phi$ satisfying
\begin{equation}
\|\phi\|_{H^2_{\mathrm{pw}}(\Omega)}
\le
C_{\mathrm{reg}}
\|z\|_{L^2(\Omega)}.
\label{e.piecewise_regularity}
\end{equation}
Here the regularity constant $C_{\mathrm{reg}}>0$ depends on the domain and interface geometry, the coefficient regularity, and the ellipticity constants, but is independent of $z$.
\end{lemma}
With this regularity estimate in hand, we now introduce the residual quantities that measure the extent to which a piecewise approximation fails to satisfy the governing equation and the associated boundary and interface conditions.

For a piecewise $H^2$ approximation
$\widetilde u=\{\widetilde u_j\}_{j=0}^K$, define the subdomain residuals
\begin{equation}
R_j(\widetilde u)=-\nabla\cdot
\left(
\beta_j\nabla\widetilde u_j
\right)
-f_j,~j=0,\ldots,K,
\end{equation}
the boundary residual
\begin{equation}
B(\widetilde u)=\widetilde u_0-g,
\end{equation}
and the interface residuals
\begin{align}
J_k(\widetilde u)&=\widetilde u_{i_k^{+}}-\widetilde u_{i_k^{-}}-v_k,
\\
F_k(\widetilde u)&=
\left(
\beta_{i_k^{+}}\nabla\widetilde u_{i_k^{+}}-\beta_{i_k^{-}}\nabla\widetilde u_{i_k^{-}}
\right)\cdot\mathbf n_k-w_k.
\end{align}

We define the continuous physics-informed residual functional by
\begin{equation*}
\mathcal J_{\mathrm c}(\widetilde u)=\sum_{j=0}^{K}
\|R_j(\widetilde u)\|_{L^2(\Omega_j)}^2+\|B(\widetilde u)\|_{L^2(\partial\Omega)}^2
+\sum_{k=1}^{K}\|J_k(\widetilde u)\|_{L^2(\Gamma_k)}^2+\sum_{k=1}^{K}
\|F_k(\widetilde u)\|_{L^2(\Gamma_k)}^2.
\end{equation*}

\begin{theorem}[Residual-to-solution stability]
\label{thm:residual_stability}
Let $u=\{u_j\}_{j=0}^K$ be the exact solution of the elliptic interface problem, and let $\widetilde u=\{\widetilde u_j\}_{j=0}^K$ be any piecewise $H^2$ approximation. Under Assumption~\ref{ass:regularity}, there exists a constant $C>0$, independent of $\widetilde u$, such that
\begin{equation}
\|\widetilde u-u\|_{L^2(\Omega)}\le
C\,\mathcal J_{\mathrm c}(\widetilde u)^{1/2}.
\label{e.residual_stability}
\end{equation}
The constant $C$ depends on the ellipticity constants, the geometry of the domain and interfaces, and the regularity constant
$C_{\mathrm{reg}}$, but is independent of the particular approximation
$\widetilde u$.
\end{theorem}

\begin{corollary}
\label{cor:dsan_stability}
Let $u_\theta=\{u_{\theta,j}\}_{j=0}^K$ be the DSAN approximation. Under the assumptions of Theorem~\ref{thm:residual_stability},
\begin{equation}
\|u_\theta-u\|_{L^2(\Omega)}
\le
C\,\mathcal J_{\mathrm c}(u_\theta)^{1/2}.
\end{equation}
Consequently, convergence of the continuous physics-informed residual to zero implies convergence of the DSAN approximation to the exact interface solution in  $L^2$ norm.
\end{corollary}


Theorem~\ref{thm:residual_stability} is a deterministic stability result formulated in terms of the continuous residual functional $\mathcal J_{\mathrm c}$. It shows that controlling the PDE, boundary, and interface residuals in their continuous norms is sufficient to control the solution error. In practice, however, DSAN is trained by minimizing residuals evaluated at only finitely many randomly sampled collocation points. Therefore, the stability estimate alone does not establish an error bound for the trained network. To bridge this gap, the next subsection relates the empirical training loss to a continuous population loss and subsequently accounts for approximation and optimization errors.

\subsection{Generalization and Approximation Error}


The physics-informed objective introduced in Section~\ref{ss.network} is evaluated on finitely many collocation points and therefore serves as an empirical approximation to a corresponding continuous loss. To make this connection precise, we first define the normalized population loss
\begin{align}
\mathcal L_{\mathrm c}(\theta)
={}&
\frac{\lambda_{\mathrm{PDE}}}{K+1}
\sum_{j=0}^{K}
\mathbb E_{\mathbf X_j\sim\mathrm{Unif}(\Omega_j)}
\left[
|R_j(u_\theta)(\mathbf X_j)|^2
\right]\\
&+
\lambda_{\mathrm{BC}}
\mathbb E_{\mathbf X_b\sim\mathrm{Unif}(\partial\Omega)}
\left[
|B(u_\theta)(\mathbf X_b)|^2
\right]
\nonumber\\
&+
\frac{\lambda_{\mathrm{jump}}}{K}
\sum_{k=1}^{K}
\mathbb E_{\mathbf X_{\Gamma_k}\sim\mathrm{Unif}(\Gamma_k)}
\left[
|J_k(u_\theta)(\mathbf X_{\Gamma_k})|^2
\right]
\nonumber\\
&+
\frac{\lambda_{\mathrm{flux}}}{K}
\sum_{k=1}^{K}
\mathbb E_{\mathbf X_{\Gamma_k}\sim\mathrm{Unif}(\Gamma_k)}
\left[
|F_k(u_\theta)(\mathbf X_{\Gamma_k})|^2
\right]\nonumber
.
\label{e.population_loss}
\end{align}
Here, the uniform distributions on $\partial\Omega$ and $\Gamma_k$ are understood with respect to the corresponding surface measures.

The population loss $\mathcal L_{\mathrm c}$ is the normalized and weighted
counterpart of the continuous residual functional $\mathcal J_{\mathrm c}$
introduced in Section~\ref{ss.residual}. For example,
\begin{equation*}
\mathbb E_{\mathbf X_j\sim\mathrm{Unif}(\Omega_j)}
\left[
|R_j(u_\theta)(\mathbf X_j)|^2
\right]
=
\frac{1}{|\Omega_j|}
\|R_j(u_\theta)\|_{L^2(\Omega_j)}^2,
\end{equation*}
and analogous identities hold for the boundary and interface residuals.
Hence, for fixed positive loss weights and fixed domain and interface
measures, there exist constants $c_1,c_2>0$ such that
\(
c_1\,\mathcal J_{\mathrm c}(u_\theta)
\le
\mathcal L_{\mathrm c}(\theta)
\le
c_2\,\mathcal J_{\mathrm c}(u_\theta).
\label{e:Jc_Lc_equivalence}
\)
Thus, controlling the population loss is equivalent, up to fixed multiplicative constants, to controlling the continuous residual functional appearing in the stability estimate.

We denote the finite-sample physics-informed training objective defined in Section~\ref{ss.network} by $\widehat{\mathcal L}_N(\theta)$. Explicitly,
\begin{align}
\widehat{\mathcal L}_N(\theta)
=&
\frac{\lambda_{\mathrm{PDE}}}{K+1}
\sum_{j=0}^{K}
\frac{1}{N_j}
\sum_{i=1}^{N_j}
|R_j(u_\theta)(\mathbf x_{j,i})|^2+
\frac{\lambda_{\mathrm{BC}}}{N_b}
\sum_{i=1}^{N_b}
|B(u_\theta)(\mathbf x_{b,i})|^2
\nonumber\\
&+
\frac{\lambda_{\mathrm{jump}}}{K}
\sum_{k=1}^{K}
\frac{1}{N_{\Gamma_k}}
\sum_{i=1}^{N_{\Gamma_k}}
|J_k(u_\theta)(\mathbf x^\Gamma_{k,i})|^2
+
\frac{\lambda_{\mathrm{flux}}}{K}
\sum_{k=1}^{K}
\frac{1}{N_{\Gamma_k}}
\sum_{i=1}^{N_{\Gamma_k}}
|F_k(u_\theta)(\mathbf x^\Gamma_{k,i})|^2 .
\label{e.empirical_loss_analysis}
\end{align}
Thus, $\widehat{\mathcal L}_N$ is precisely the empirical training loss, whereas $\mathcal L_{\mathrm c}$ is its population counterpart. The next step is to quantify, uniformly over the admissible parameter set $\Theta$, the discrepancy between these two quantities.

For each subdomain $\Omega_j$, define
\(
\mathcal G_{\mathrm{PDE},j}=
\left\{\mathbf x\mapsto|R_j(u_\theta)(\mathbf x)|^2:\theta\in\Theta
\right\},
\)
and define analogously $\mathcal G_{\mathrm{BC}}$, $\mathcal G_{\mathrm{jump},k}$, and $\mathcal G_{\mathrm{flux},k}$.
Let $\mathfrak R_N(\mathcal G)$ denote the Rademacher complexity of a function class $\mathcal G$ based on $N$ independent samples. 
To relate the empirical training loss to its population counterpart, we make the following assumptions on the sampling procedure and residual classes.

\begin{assumption}[Sampling and residual classes]
\label{ass:sampling}
For each $j=0,\ldots,K$, the collocation points $\{\mathbf x_{j,i}\}_{i=1}^{N_j}$ are sampled independently and uniformly
from $\Omega_j$. The boundary points are sampled independently and uniformly from $\partial\Omega$. For each $k=1,\ldots,K$, the interface points $\{\mathbf x^\Gamma_{k,i}\}_{i=1}^{N_{\Gamma_k}}$ are sampled independently and uniformly from $\Gamma_k$. All collocation points are generated once before training and remain fixed throughout the optimization.
The admissible parameter set $\Theta$ is compact. Moreover, for each residual class introduced above, there exists a finite constant $M_s$ such that
\(
0\le g(\mathbf x)\le M_s,
\)
for all $g\in\mathcal G_s$, $s\in\{PDE,~BC,~jump,~flux\}$, and all points in the corresponding sampling domain.
\end{assumption}

Then we have the following estimation between the empirical loss and the population loss.

\begin{theorem}[Uniform empirical-to-population bound]
\label{thm:uniform_loss}
Under Assumption~\ref{ass:sampling}, for any $\delta\in(0,1)$, with probability at least $1-\delta$ over the initial draw of the collocation points,
\begin{equation}
\sup_{\theta\in\Theta}
\left|\mathcal L_{\mathrm c}(\theta)-\widehat{\mathcal L}_{N}(\theta)\right|\le\Delta_N(\delta),
\label{e.empirical_population}
\end{equation}
where
\(
\Delta_N(\delta)=
\Delta_{\mathrm{PDE}}+
\Delta_{\mathrm{BC}}+
\Delta_{\mathrm{jump}}+
\Delta_{\mathrm{flux}},
\)
with
\begin{align}
\Delta_{\mathrm{PDE}}={}&\frac{\lambda_{\mathrm{PDE}}}{K+1}
\sum_{j=0}^{K}
\left[2\mathfrak R_{N_j}(\mathcal G_{\mathrm{PDE},j})+M_{\mathrm{PDE},j}\sqrt{\frac{\log(2Q/\delta)}{2N_j}}\right],\nonumber\\
\Delta_{\mathrm{BC}}={}&\lambda_{\mathrm{BC}}
\left[
2\mathfrak R_{N_b}(\mathcal G_{\mathrm{BC}})+M_{\mathrm{BC}}\sqrt{\frac{\log(2Q/\delta)}{2N_b}}
\right],
\nonumber\\
\Delta_{\mathrm{jump}}={}&\frac{\lambda_{\mathrm{jump}}}{K}\sum_{k=1}^{K}\left[2\mathfrak R_{N_{\Gamma_k}}
(\mathcal G_{\mathrm{jump},k})+M_{\mathrm{jump},k}\sqrt{\frac{\log(2Q/\delta)}{2N_{\Gamma_k}}}
\right],
\nonumber\\
\Delta_{\mathrm{flux}}
={}&
\frac{\lambda_{\mathrm{flux}}}{K}
\sum_{k=1}^{K}
\left[
2\mathfrak R_{N_{\Gamma_k}}
(\mathcal G_{\mathrm{flux},k})
+
M_{\mathrm{flux},k}
\sqrt{
\frac{\log(2Q/\delta)}
{2N_{\Gamma_k}}
}
\right].\nonumber
\end{align}
Here $M_s$ denotes a uniform upper bound for the corresponding squared residual class and  
\(Q=3K+2.\)
\end{theorem}

In particular, if $\mathfrak R_{N_s}(\mathcal G_s)\to0$ for every residual class as $N_s\to\infty$, then $\Delta_N(\delta)\to0$, yielding uniform consistency of the empirical loss. However, Theorem~\ref{thm:uniform_loss} controls only the error introduced by replacing continuous expectations with finite-sample averages. To characterize the network actually returned by training, we must additionally account for the fact that the optimization algorithm may not attain the global minimum of the empirical loss. This motivates the following optimization-error formulation.

\begin{corollary}[Population loss of the trained network]
\label{cor:population_loss}
Suppose that the trained parameter $\widehat\theta$ satisfies
\begin{equation}
\widehat{\mathcal L}_{N}(\widehat\theta)
\le
\inf_{\theta\in\Theta}
\widehat{\mathcal L}_{N}(\theta)
+
\varepsilon_{\mathrm{opt}},
\label{e.optimization_error}
\end{equation}
where $\varepsilon_{\mathrm{opt}}$ denotes the optimization error, measuring the gap between the empirical loss attained by the training
algorithm and the minimum empirical loss over the admissible DSAN hypothesis class. Then, on the event of Theorem~\ref{thm:uniform_loss},
\begin{equation}
\mathcal L_{\mathrm c}(\widehat\theta)
\le
\inf_{\theta\in\Theta}\mathcal L_{\mathrm c}(\theta)+\varepsilon_{\mathrm{opt}}+2\Delta_N(\delta).
\label{e.population_bound}
\end{equation}
\end{corollary}
Thus, the population loss of the trained network is controlled by the best population loss achievable within the DSAN hypothesis class, together with the optimization error and the finite-sample generalization error.

It remains to quantify the first term on the right-hand side of \eqref{e.population_bound}, namely the smallest population loss attainable within the DSAN hypothesis class. This term reflects the approximation capability of the network architecture. We therefore introduce the piecewise $H^2$ approximation error
\(
\varepsilon_{\mathrm{app}}=
\inf_{\theta\in\Theta}
\|u_\theta-u\|_{H^2_{\mathrm{pw}}(\Omega)}.
\label{e.dsan_app_error}
\)

\begin{lemma}[Approximation error controls the population residual]
\label{lem:approximation_residual}
Under Assumption~\ref{ass:regularity}, there exists a constant $C_{\mathrm{app}}>0$ such that, for every $u_\theta\in H^2_{\mathrm{pw}}(\Omega)$,
\begin{equation}
\mathcal L_{\mathrm c}(\theta)
\le
C_{\mathrm{app}}
\|u_\theta-u\|_{H^2_{\mathrm{pw}}(\Omega)}^2.
\label{e.app_residual_bound}
\end{equation}
Consequently,
\begin{equation}
\inf_{\theta\in\Theta}
\mathcal L_{\mathrm c}(\theta)
\le
C_{\mathrm{app}}
\varepsilon_{\mathrm{app}}^2.
\end{equation}
\end{lemma}
Lemma~\ref{lem:approximation_residual} therefore converts approximation in the piecewise $H^2$ norm into control of the continuous physics-informed loss. Together with Corollary~\ref{cor:population_loss}, it yields
\[
\mathcal L_{\mathrm c}(\widehat\theta)\le
C_{\mathrm{app}}\varepsilon_{\mathrm{app}}^2+\varepsilon_{\mathrm{opt}}+2\Delta_N(\delta).
\]

The approximation error $\varepsilon_{\mathrm{app}}$ above is defined
abstractly for the complete DSAN hypothesis class and does not yet expose the
role of the shared backbone and the subdomain-specific adapters. To connect
the error analysis more directly to the proposed architecture, we next
decompose both the exact solution and the DSAN approximation into shared and
domain-specific components.
\begin{proposition}[Shared-plus-adaptation approximation]
\label{prop:shared_adapter}
The DSAN approximation in each subdomain admits the decomposition
\( 
u_{\theta,j}=
u_{\theta}^{\mathrm{sh}}+
u_{\theta,j}^{\mathrm{ad}},~j=0,\ldots,K,
\)
as introduced in \eqref{e.dsan_decomposition}.

Let $u^{\mathrm{sh}}$ be a common reference component defined on $\Omega$, and decompose the exact solution in each subdomain as
\( 
u_j=u^{\mathrm{sh}}|_{\Omega_j}+u_j^{\mathrm{ad}},~
u_j^{\mathrm{ad}}=u_j-u^{\mathrm{sh}}|_{\Omega_j}.
\label{e.exact_shared_adaptation}
\)
Then
\begin{equation}
\|u-u_\theta\|_{H^2_{\mathrm{pw}}(\Omega)}^2
\le
\sum_{j=0}^{K}\left(
\|u^{\mathrm{sh}}-u_{\theta}^{\mathrm{sh}}\|_{H^2(\Omega_j)}+
\|u_j^{\mathrm{ad}}-u_{\theta,j}^{\mathrm{ad}}\|_{H^2(\Omega_j)}
\right)^2.
\label{e.shared_adapter_bound}
\end{equation}
\end{proposition}

Proposition~\ref{prop:shared_adapter} shows that the overall approximation error can be controlled by two complementary representation errors. The error in the common shared component and the errors in the subdomain-specific adaptations. This motivates the architecture-aware approximation measure $\varepsilon_{\mathrm{SA}}$.
\begin{equation}
\varepsilon_{\mathrm{SA}}=
\inf_{\theta\in\Theta}
\left[
\sum_{j=0}^{K}
\left(
\|u^{\mathrm{sh}}-u_{\theta}^{\mathrm{sh}}\|_{H^2(\Omega_j)}+
\|u_k^{\mathrm{ad}}-u_{\theta,j}^{\mathrm{ad}}\|_{H^2(\Omega_j)}
\right)^2
\right]^{1/2}.
\label{e.sa_error}
\end{equation}

It follows directly from Proposition~\ref{prop:shared_adapter} that
\(
\varepsilon_{\mathrm{app}}
\le
\varepsilon_{\mathrm{SA}}.
\label{e.app_sa_relation}\)
We now have all ingredients needed to bound the error of the trained DSAN. Theorem~\ref{thm:residual_stability} converts the continuous residual into solution error, Theorem~\ref{thm:uniform_loss} and Corollary~\ref{cor:population_loss} control the population loss from finite-sample training, and Lemma~\ref{lem:approximation_residual} bounds the best attainable population loss by the network approximation error. Proposition~ \ref{prop:shared_adapter} further provides a DSAN-specific interpretation of this approximation term. Combining these results yields the following total error estimate.

\begin{theorem}[Total error bound for the trained DSAN]
\label{thm:total_error}
Suppose Assumptions~\ref{ass:regularity} and~\ref{ass:sampling} hold.
Let $\widehat\theta$ be the parameter obtained by training and assume that
\begin{equation*}
\widehat{\mathcal L}_{N}(\widehat\theta)
\le
\inf_{\theta\in\Theta}
\widehat{\mathcal L}_{N}(\theta)+
\varepsilon_{\mathrm{opt}}.
\label{e.optimization_error2}
\end{equation*}
Then, for any $\delta\in(0,1)$, with probability at least $1-\delta$
over the initial random draw of the fixed collocation points,
\begin{equation}
\|u_{\widehat\theta}-u\|_{L^2(\Omega)}
\le
C\left[C_{\mathrm{app}}\varepsilon_{\mathrm{app}}^2+\varepsilon_{\mathrm{opt}}+2\Delta_N(\delta)
\right]^{1/2}.
\label{e.total_error_bound}
\end{equation}
In particular, using
$\varepsilon_{\mathrm{app}}\le\varepsilon_{\mathrm{SA}}$,
\begin{equation}
\|u_{\widehat\theta}-u\|_{L^2(\Omega)}
\le
C
\left[C_{\mathrm{app}}\varepsilon_{\mathrm{SA}}^2+\varepsilon_{\mathrm{opt}}+2\Delta_N(\delta)
\right]^{1/2}.
\label{e.total_error_sa}
\end{equation}
\end{theorem}

Estimate~\eqref{e.total_error_sa} makes explicit the three sources of error in the trained DSAN. The term $\varepsilon_{\mathrm{SA}}$ represents the approximation error associated with the shared-backbone and domain-specific adapter architecture, $\varepsilon_{\mathrm{opt}}$ measures the optimization error incurred during training, and $\Delta_N(\delta)$ quantifies the finite-sample generalization error. The stability constant transfers these three loss-level contributions to the physical solution error in $L^2(\Omega)$.

\subsection{Parameter Complexity Analysis} 
The preceding analysis addresses the solution accuracy of the trained DSAN. We now turn to the second theoretical aspect of the proposed architecture, namely its parameter efficiency.  Specifically, we compare the growth in the number of trainable parameters of DSAN with that of a fully decomposed domain-specific formulation, in which each subdomain is assigned an independent neural network.

Consider a fully connected neural network with input dimension $d$, hidden width $D$, $L$ hidden layers, and a scalar output.
For a multi-layer perceptron, the input layer contains \(P_{\mathrm{in}} = dD + D\) trainable parameters, including the bias terms. The transformations between successive hidden layers contain \(P_{\mathrm{hid}}=(L-1)(D^2+D)\) parameters, while the output layer contains \(P_{\mathrm{out}}=D+1.\)
Therefore, the total number of parameters in one complete network is
\begin{equation*}
P_{\mathrm{net}}=(L-1)D^2+(d+L+1)D+1.
\label{e.pnet}
\end{equation*}

A conventional domain-decomposition method assigns an independent complete network to each of the $K+1$ subdomains. Its total number of trainable parameters is therefore
\begin{equation}
P_{\mathrm{DD}} = (K+1)P_{\mathrm{net}}.
\label{e.pdd}
\end{equation}

In DSAN, the backbone and output layer are shared across all subdomains, while each subdomain introduces only a lightweight adapter. For an adapter with bottleneck dimension $r$, the down-projection and up-projection contain
\begin{equation*}
P_{\mathrm a} = (rD+r)+(Dr+D) = 2Dr+D+r
\label{e.pa}
\end{equation*}
trainable parameters. Consequently, the total number of trainable parameters in DSAN is
\begin{equation}
P_{\mathrm{DSAN}} = P_{\mathrm{net}} + (K+1)P_{\mathrm a}.
\label{e.pdsan}
\end{equation}
Thus, adding a new subdomain requires an additional complete network with $P_{\mathrm{net}}$ parameters in the conventional domain-decomposition formulation, whereas DSAN requires only an additional adapter with $P_{\mathrm a}$ parameters.

The relative parameter complexity of the two formulations is
\begin{equation}
\frac{P_{\mathrm{DSAN}}}{P_{\mathrm{DD}}} = \frac{P_{\mathrm{net}}+(K+1)P_{\mathrm a}}{(K+1)P_{\mathrm{net}}}
=\frac{1}{K+1}+\frac{P_{\mathrm a}}{P_{\mathrm{net}}}.
\label{e.parameter_ratio}
\end{equation}
Since $r\ll D$, the adapter contains substantially fewer parameters than a complete network. In particular,
\begin{equation*}
\begin{aligned}
\frac{P_{\mathrm a}}{P_{\mathrm{net}}}&=
\frac{2Dr+D+r}{(L-1)D^2+(d+L+1)D+1}
<\frac{(2r+2)D}{(L-1)D^2}=\frac{2r+2}{(L-1)D},
\end{aligned}
\label{e.adapter_ratio}
\end{equation*}
where the inequality holds for $D>r$.

Hence, for fixed bottleneck dimension $r$ and network depth $L$, the relative cost of a subdomain-specific adapter decreases as the backbone width $D$ increases. For example, when $r=4$ and $L\geq3$,
\(
\frac{P_{\mathrm a}}{P_{\mathrm{net}}}<\frac{5}{D}.
\)
Therefore, the parameter overhead introduced by each additional subdomain is small compared with replicating a complete network.

The different scaling behavior becomes particularly clear as the number of subdomains increases. From \eqref{e.pdd} and \eqref{e.pdsan},
\begin{equation}
P_{\mathrm{DD}}=\mathcal O\!\left((K+1)LD^2\right),
\end{equation}
whereas
\begin{equation}
P_{\mathrm{DSAN}}=
\mathcal O\!\left(LD^2+(K+1)Dr\right).
\end{equation}
Thus, the incremental parameter cost per additional subdomain is $\mathcal O(LD^2)$ for the fully decomposed formulation but only $\mathcal O(Dr)$ for DSAN.
Moreover,
\begin{equation}
\lim_{K\rightarrow\infty}
\frac{P_{\mathrm{DSAN}}}{P_{\mathrm{DD}}}=
\frac{P_{\mathrm a}}{P_{\mathrm{net}}}.
\label{e.parameter_limit}
\end{equation}
This asymptotic relation characterizes the parameter scaling when the number of subdomains becomes large. As $K$ increases, the fixed cost of the shared backbone is amortized over an increasing number of subdomains, and the relative parameter ratio approaches $\frac{P_{\mathrm a}}{P_{\mathrm{net}}}$. Since $P_{\mathrm a}\ll P_{\mathrm{net}}$ for $r\ll D$, DSAN retains a substantially smaller parameter count than a fully decomposed formulation for problems with many subdomains.
\section{Experiments}
\label{s.exps}

\subsection{Experimental Setup}
We consider four benchmark problems covering two- and three-dimensional domains, isotropic and anisotropic coefficients, and configurations with 1 to 27 interfaces. DSAN is compared with XPINN~\cite{jagtap2020extended}, which assigns an independent network to each subdomain, and CEPINN~\cite{hu2026categorical}, which employs a single network with categorical embedding. For each benchmark, all methods use the same collocation, boundary, and interface points, as well as the same loss functions and corresponding loss weights, to ensure a consistent comparison.

We evaluate the solution accuracy using the mean squared error (MSE) over the computational domain, together with the MSEs of the solution jump and flux jump across the interfaces. To assess scalability, we further report the number of trainable parameters for each method.

All neural models are trained using a two-stage optimization procedure. We first employ the Adam optimizer with a learning rate of $0.001$ for 20000 iterations, followed by 5000 iterations of L-BFGS optimization.
The hyperbolic tangent (tanh) function is used as the activation function throughout all experiments. For DSAN, the adapter bottleneck dimension $r$ is fixed at 4. For CEPINN, the embedding dimension is set to 2.
The network hyperparameters used for each benchmark are summarized in Table~\ref{T.param}. Unless otherwise specified, the compared methods use the same network hyperparameters to provide a controlled comparison. All experiments are conducted on a workstation equipped with an NVIDIA GeForce RTX 4090 GPU.

\begin{table}[ht]
\centering
\caption{Summary of the experimental setup. \#Subdomain Points denotes the total number of collocation points sampled in the subdomains; \#Interface Points denotes the total number of points sampled on all interfaces; and \#Boundary Points denotes the total number of points sampled on the domain boundary. Depth is the number of layers in an MLP, while width is the number of neurons in each hidden layer.}
\begin{adjustbox}{width=\textwidth}
\begin{tabular}{lcccccc}
\toprule[2pt]
Problem & PDE Dimension & Depth & Width & \#Subdomain Points &
\#Interface Points & \#Boundary Points \\
\midrule[1pt]
Example 1 & 3 & 4 & 64  & 10000 & 1000  & 6000 \\
Example 2 & 2 & 4 & 64  & 12000 & 2000  & 1600 \\
Example 3 & 2 & 4 & 64  & 9500  & 1800  & 2000 \\
Example 4 & 3 & 4 & 128 & 33500 & 13824 & 6000 \\
\bottomrule[2pt]
\end{tabular}
\end{adjustbox}
\label{T.param}
\end{table}

\subsection{Example 1: 3D ellipsoidal interface problem}

\begin{figure}[!ht]
    \centering
    \includegraphics[width=0.5\linewidth]{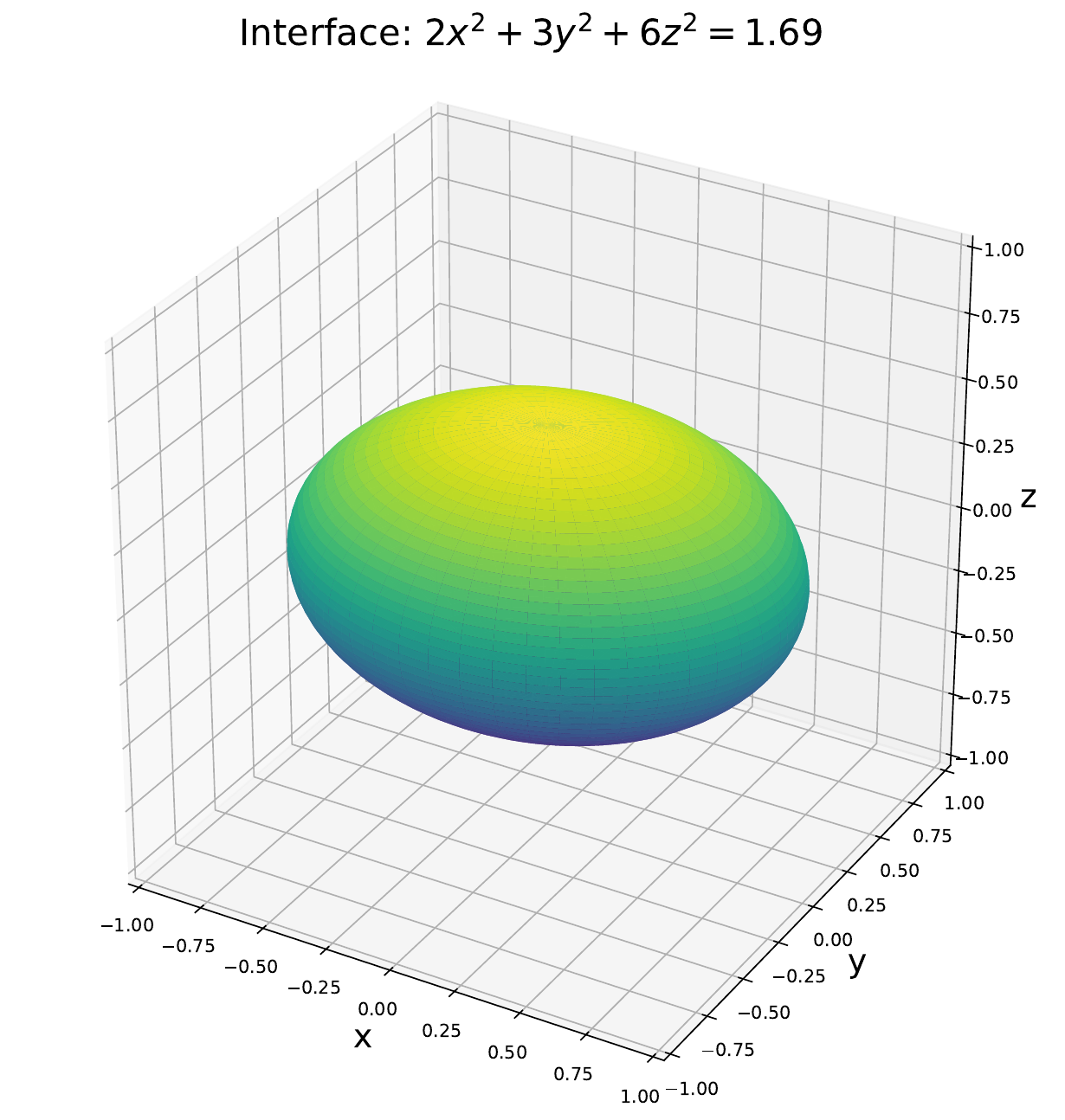}
    \caption{3D ellipsoidal interface visualization.}
    \label{f.case1interface}
\end{figure}

We first consider a three-dimensional interface problem defined on $\Omega=[-1,1]^3$. The domain is separated into two subdomains by an ellipsoidal interface
\begin{equation*}
\Gamma=\left\{(x,y,z)\in\Omega:2x^2+3y^2+6z^2=1.69\right\},
\end{equation*}
as illustrated in Fig.~\ref{f.case1interface}. 
A piecewise-constant scalar coefficient $\beta(x)$ is defined as follows 
\begin{equation*}
\beta(\mathbf{x})=
\begin{cases}
2, & 2x^2+3y^2+6z^2\leq 1.69, \\
60, & 2x^2+3y^2+6z^2>1.69.
\end{cases}
\end{equation*}
The corresponding exact solution is defined piecewise as
\begin{equation}
u(\mathbf{x})=
\begin{cases}
x^2+y^2+z^2, & 2x^2+3y^2+6z^2\leq 1.69, \\
x+y+z, & 2x^2+3y^2+6z^2>1.69.
\end{cases}
\label{e:ex1_exact}
\end{equation}
and the boundary condition is $g(\mathbf{x})=x+y+z$.
The solution therefore exhibits different functional behaviors on the two sides of the interface and is discontinuous across $\Gamma$. 

Table~\ref{tab:ex1} reports the quantitative results of the three methods. All three formulations accurately approximate the piecewise solution and achieve domain MSEs of the same order of magnitude, with DSAN obtaining the lowest error, XPINN yielding a comparable value, and CEPINN producing a slightly larger error. 
The jump MSEs exhibit a similar trend, with DSAN again attaining the lowest error. The difference is more pronounced in the flux MSE, where DSAN achieves approximately half the values obtained by CEPINN and XPINN. DSAN achieves these improvements with only slightly more trainable parameters than CEPINN and substantially fewer than XPINN. This result indicates that the domain-specific adapters can effectively represent the different gradient behaviors on the two sides of the interface while maintaining a shared backbone.

\begin{figure}[!ht]
    \centering
    \includegraphics[width=0.9\linewidth]{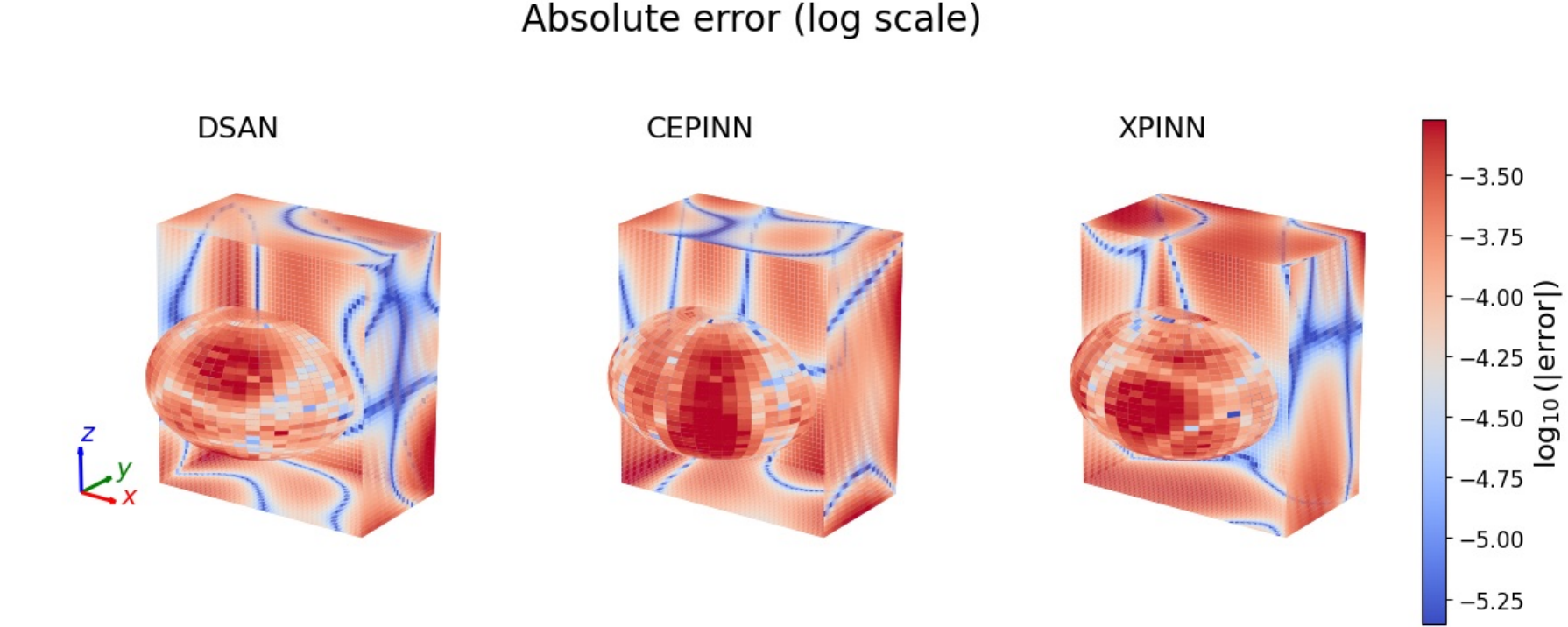}
    \caption{Cross-sectional comparison of absolute solution errors for the 3D ellipsoidal interface problem.}
    \label{f.case1comparison}
\end{figure}

The absolute-error distributions are shown in Fig.~\ref{f.case1comparison}. The three methods maintain small errors throughout
the computational domain, although the spatial distributions of the errors differ. Together with the quantitative results, these observations show that DSAN achieves accuracy comparable to the fully decomposed XPINN while requiring fewer trainable parameters. Thus, for this single-interface problem, DSAN retains the approximation capability without requiring a complete independent network for each subdomain.

\begin{table}[!ht]
\centering
\caption{Accuracy and trainable parameter comparison for the 3D ellipsoidal interface problem. The best result for each measurement is highlighted in \textbf{bold}.}
\label{tab:ex1}
\begin{tabular}{lcccc}
\toprule
Method & Parameters & MSE & Jump MSE & Flux MSE\\
\midrule
DSAN & 11359 & \textbf{3.55$\times$10$^{-8}$} & \textbf{1.85$\times$10$^{-7}$} & \textbf{8.57$\times$10$^{-5}$}\\
CEPINN  & \textbf{10551} & $4.60{\times}10^{-8}$ & $3.18{\times}10^{-7}$ & 1.95$\times$10$^{-4}$\\
XPINN   & 20902 & 3.78$\times$10$^{-8}$ & 2.96$\times$10$^{-7}$ & $1.93{\times}10^{-4}$\\
\bottomrule
\end{tabular}
\end{table}

\subsection{Example 2: Multiple Concentric Interfaces}

We consider $\Omega=[-1,1]^2$ with five concentric circular interfaces
$$\Gamma_k =\left\{(x,y)\in\Omega: x^2+y^2=R_k^2 \right\},~k=1,\ldots,5,$$
with $(R_1,R_2,R_3,R_4,R_5)=(0.80,0.65,0.50,0.35,0.20).$
The interfaces partition $\Omega$ into six subdomains
$\{\Omega_k\}_{k=0}^5$, numbered from the exterior to the interior, so that
$\Gamma_k$ separates $\Omega_{k-1}$ and $\Omega_k$. We take
$i_k^+=k$, $i_k^-=k-1$, with the unit normal
$\mathbf n_k=\frac{(x, y)}{R_k}$ pointing from $\Omega_k$ to $\Omega_{k-1}$.

The diffusion coefficient is piecewise constant,
\[
\beta|_{\Omega_k}=\beta_k,~
(\beta_0,\ldots,\beta_5)=(4,20,2,10,5,1),
\]
and the exact solution is
\begin{equation}
u_k(\mathbf x)=\frac{x^2+y^2}{\beta_k}+b_k,~
\mathbf x\in\Omega_k.
\label{e.case2_solution}
\end{equation}
We set $b_5=0$ and determine the remaining constants recursively by
$$b_{k-1}=b_k+R_k^2\left(\frac{1}{\beta_k}-\frac{1}{\beta_{k-1}}\right)-v_k,~k=5,\ldots,1,$$
where $(v_1,\ldots,v_5)=(-0.6,\,0.2,\,-0.4,\,0.3,\,-0.5).$
With this construction,
\begin{equation}
[u]_{\Gamma_k}=v_k,~~[\beta\nabla u\cdot\mathbf n_k]_{\Gamma_k}=0,~k=1,\ldots,5,
\end{equation}
and the source term is $f=-4$ in all subdomains.

\begin{figure}[!ht]
    \centering
    \includegraphics[width=0.75\linewidth]{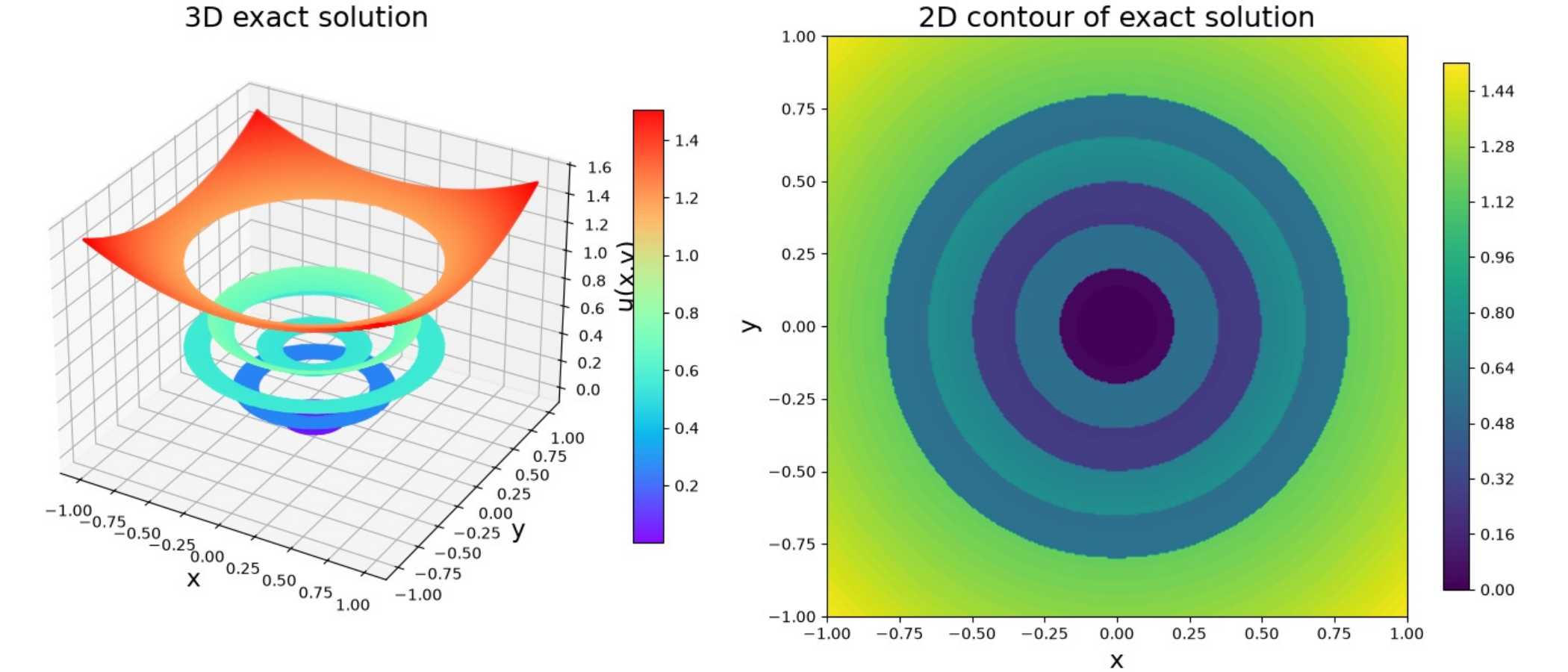}
    \caption{Exact solution and its contour representation for the 2D problem with five concentric interfaces.}
    \label{f.case2interface}
\end{figure}

This example therefore represents a chain-type interface configuration with
discontinuous solution values and continuous normal flux across all five
interfaces. Figure~\ref{f.case2interface} shows the exact solution and its
contour representation.

\begin{table}[!ht]
\centering
\caption{Accuracy and model complexity comparison for the 2D problem with five concentric interfaces. The best result for each measurement is highlighted in \textbf{bold}.}
\label{tab:ex2}
\begin{tabular}{lcccc}
\toprule
Method & Parameters & MSE & Jump MSE & Flux MSE\\
\midrule
DSAN   & 13125 & \textbf{5.97}$\times$\textbf{10}$^{-10}$ & $4.77{\times}10^{-10}$ & $6.73{\times}10^{-8}$\\
CEPINN & \textbf{10513} & $4.68{\times}10^{-9}$  & $1.24{\times}10^{-9}$  & $1.23{\times}10^{-7}$\\
XPINN  & 62406 & $1.97{\times}10^{-9}$  & \textbf{8.58}$\times$\textbf{10}$^{-11}$ & \textbf{6.64}$\times$\textbf{10}$^{-8}$\\
\bottomrule
\end{tabular}
\end{table}

A total of 2000, 400, and 400 points are sampled in each subdomain, on each interface, and on each boundary face, respectively.
Table~\ref{tab:ex2} compares the prediction accuracy and number of trainable parameters of the three methods. We can see DSAN achieves the lowest domain MSE. Thus, despite using substantially fewer parameters than XPINN, DSAN provides higher overall solution accuracy for this multi-interface problem. In particular, DSAN uses 13125 trainable parameters, compared with 62406 for XPINN, corresponding to approximately 21\% of the parameters required by the fully decomposed formulation. For the interface conditions, XPINN achieves the lowest jump MSE, although all three methods accurately reproduce the prescribed solution discontinuities. For the flux condition, DSAN and XPINN achieve comparable errors, with both outperforming CEPINN. These results indicate that DSAN can accurately represent both the discontinuous solution values and continuous normal flux across multiple interfaces without requiring independent full networks for the individual subdomains.

\begin{figure}[!ht]
    \centering
    \includegraphics[width=\linewidth]{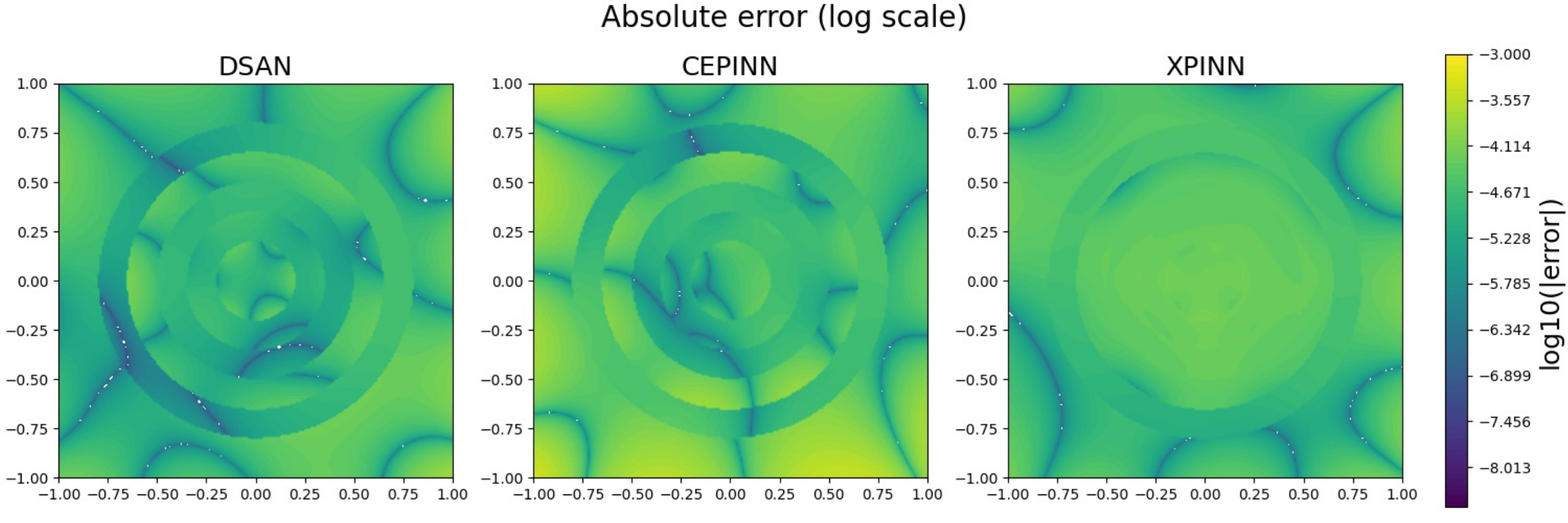}
    \caption{Comparison of the $\log_{10}$ absolute solution errors for the 2D problem with five concentric interfaces. Darker colors correspond to smaller absolute errors and therefore higher prediction accuracy.}
    \label{f.case2comparison}
\end{figure}

Fig.~\ref{f.case2comparison} further compares the spatial distributions of the absolute solution errors on a logarithmic scale. Consistent with the quantitative results, all three methods maintain small errors over most of the computational domain, although their spatial error distributions differ. DSAN exhibits a generally lower error level across large portions of the domain, whereas CEPINN shows more pronounced high-error regions. XPINN also maintains low errors over much of the domain, consistent with its competitive quantitative accuracy. Overall, the spatial error distributions further demonstrate that DSAN accurately approximates the multi-interface solution while requiring substantially fewer parameters than XPINN.

\subsection{Example 3: Anisotropic Irregular Interfaces}
\begin{figure}[!ht]
    \centering
    \includegraphics[width=\linewidth]{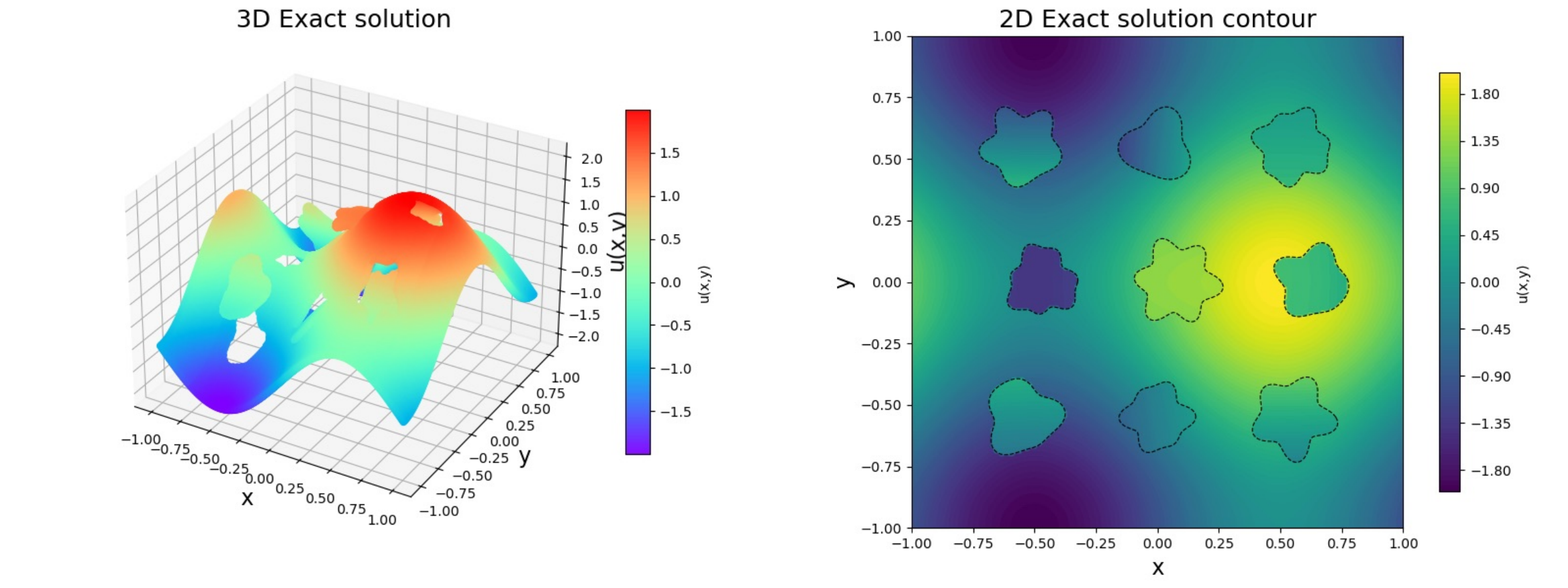}
    \caption{Exact solution for the two-dimensional anisotropic problem with nine irregular interfaces.}
    \label{f.case3interface}
\end{figure}

We next consider a 2D anisotropic interface problem with nine irregular interfaces. The domain is $\Omega=[-1,1]^2$, containing nine non-overlapping star-shaped interfaces. The interfaces divide the domain into a background region $\Omega_0$ and nine inclusion regions $\{\Omega_k\}_{k=1}^{9}$ (see Fig.~\ref{f.case3interface}).

The diffusion tensor in Eq.~\eqref{e.interface} is piecewise constant and defined as
\begin{equation*}
\mathbf \beta(\mathbf x) =
\begin{cases}
\mathbf I, & \mathbf x\in\Omega_0,\\
\mathbf A_k, & \mathbf x\in\Omega_k,
\quad k=1,\ldots,9,
\end{cases}
\end{equation*}
where $\mathbf A_k=\operatorname{diag}(a_k,b_k)$.
The diagonal entries are given by
\begin{align*}
(a_1,\ldots,a_9) = (5,1,6,0.8,8,1.2,4,0.6,10), 
(b_1,\ldots,b_9) = (1,5,0.8,6,1.2,8,0.6,4,2).
\end{align*}

The exact solution is constructed as
\begin{equation}
u(\mathbf x)=
\begin{cases}
\sin(\pi x)+\cos(\pi y), & \mathbf x\in\Omega_0, \\
\dfrac{\sin(\pi x)}{a_k}+\dfrac{\cos(\pi y)}{b_k}+c_k,& \mathbf x\in\Omega_k, ~ k=1,\ldots,9,
\end{cases}
\end{equation}
with $(c_1,\ldots,c_9)=(0.35,-0.25,0.45,-0.30,0.50,-0.20, 0.28,-0.40,0.38)$.
Accordingly, the prescribed solution jump on $\Gamma_k$ is
\begin{equation}
v_k(\mathbf x)=\left(\frac{1}{a_k}-1\right)\sin(\pi x)+\left(\frac{1}{b_k}-1\right)\cos(\pi y)+c_k.
\end{equation}
Thus, unlike the constant jumps considered in the Example 2, the solution jump varies along each irregular interface.
Moreover, the exact solution satisfies
\begin{equation}
\mathbf A_k\nabla u_k=
\mathbf A_0\nabla u_0=
\begin{pmatrix}
\pi\cos(\pi x)\\
-\pi\sin(\pi y)
\end{pmatrix},
\end{equation}
and hence the anisotropic normal flux is continuous across every interface. In contrast, the solution generally remains discontinuous across the interfaces.
Fig.~\ref{f.case3interface} illustrates the exact solution for this problem. The combination of irregular interface geometries and anisotropic diffusion coefficients introduces more complex subdomain-dependent solution behavior.

We sample 5000 points in $\Omega_0$ and 500 points in each $\Omega_k$. In addition, 2000 points are sampled on the boundary, and 200 points are sampled on each interface. Table~\ref{tab:ex3} summarizes the accuracy and model parameters of the three methods. DSAN achieves the lowest overall solution MSE, while XPINN obtains the lowest jump and flux errors, with DSAN showing slightly higher errors for these two interface conditions. In contrast, CEPINN produces larger errors across all three metrics. These results suggest that the domain-specific representations employed by DSAN and XPINN provide stronger approximation capability for this challenging anisotropic interface problem. In terms of model complexity, DSAN requires substantially fewer trainable parameters than XPINN while maintaining accurate domain and interface approximations, demonstrating that such domain-specific representation can be achieved without assigning a separate full network to each subdomain.

\begin{figure}[!ht]
    \centering
    \includegraphics[width=\linewidth]{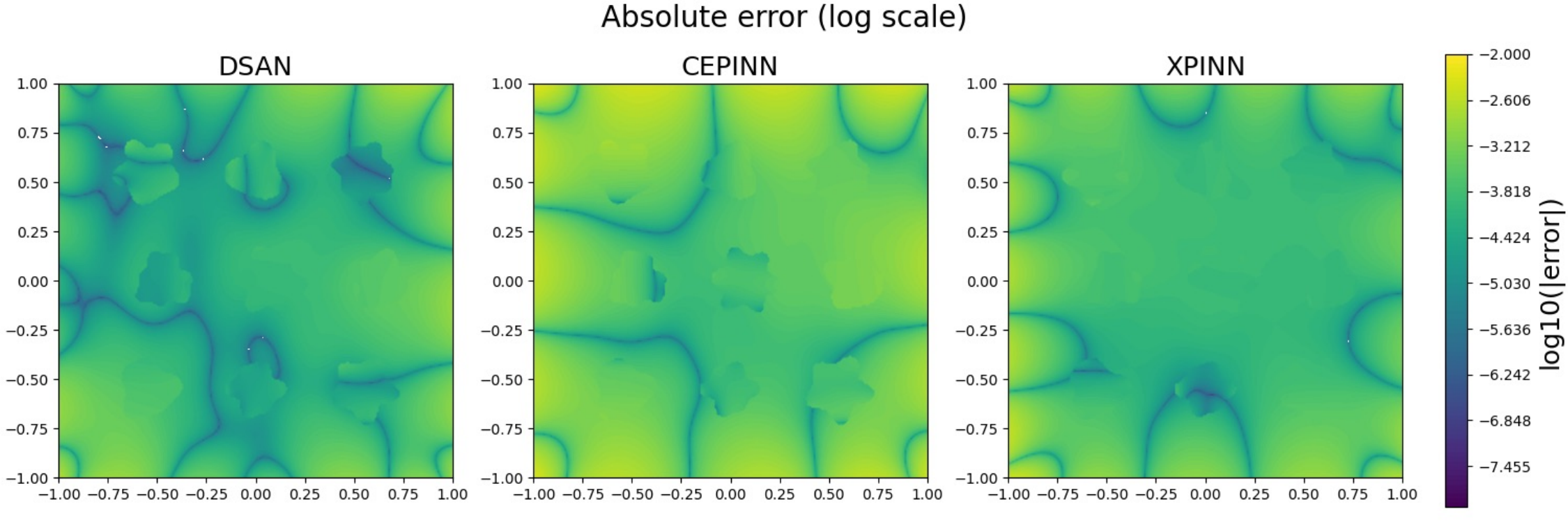}
    \caption{Comparison of absolute solution errors for the 2D anisotropic problem. Darker colors correspond to smaller absolute errors and therefore higher prediction accuracy.}
    \label{f.case3comparison}
\end{figure}

Fig.~\ref{f.case3comparison} presents the spatial distributions of the absolute solution errors. The results are consistent with the quantitative comparison, with DSAN maintaining low errors across the computational domain. This example further demonstrates the ability of DSAN to handle multiple irregular interfaces and anisotropic diffusion.

\begin{table}[!ht]
\centering
\caption{Accuracy and model complexity comparison for the 2D anisotropic problem. The best result for each measurement is highlighted in \textbf{bold}.}
\label{tab:ex3}
\begin{tabular}{lcccc}
\toprule
Method & Parameters & MSE & Jump MSE & Flux MSE\\
\midrule
DSAN   & 22697  & \textbf{6.35}$\times$\textbf{10}$^{-8}$ & $4.62{\times}10^{-9}$ & $1.99{\times}10^{-6}$\\
CEPINN & \textbf{17045}  & $6.88{\times}10^{-7}$ & $1.95{\times}10^{-8}$ & $5.87{\times}10^{-6}$\\
XPINN  & 168970 & $1.27{\times}10^{-7}$ & \textbf{1.21}$\times$\textbf{10}$^{-9}$ & \textbf{9.11}$\times$\textbf{10}$^{-7}$\\
\bottomrule
\end{tabular}
\end{table}

\subsection{Example 4: Twenty-Seven Complex Interfaces in 3D}
\begin{figure}[!ht]
    \centering
    \includegraphics[width=0.6\linewidth]{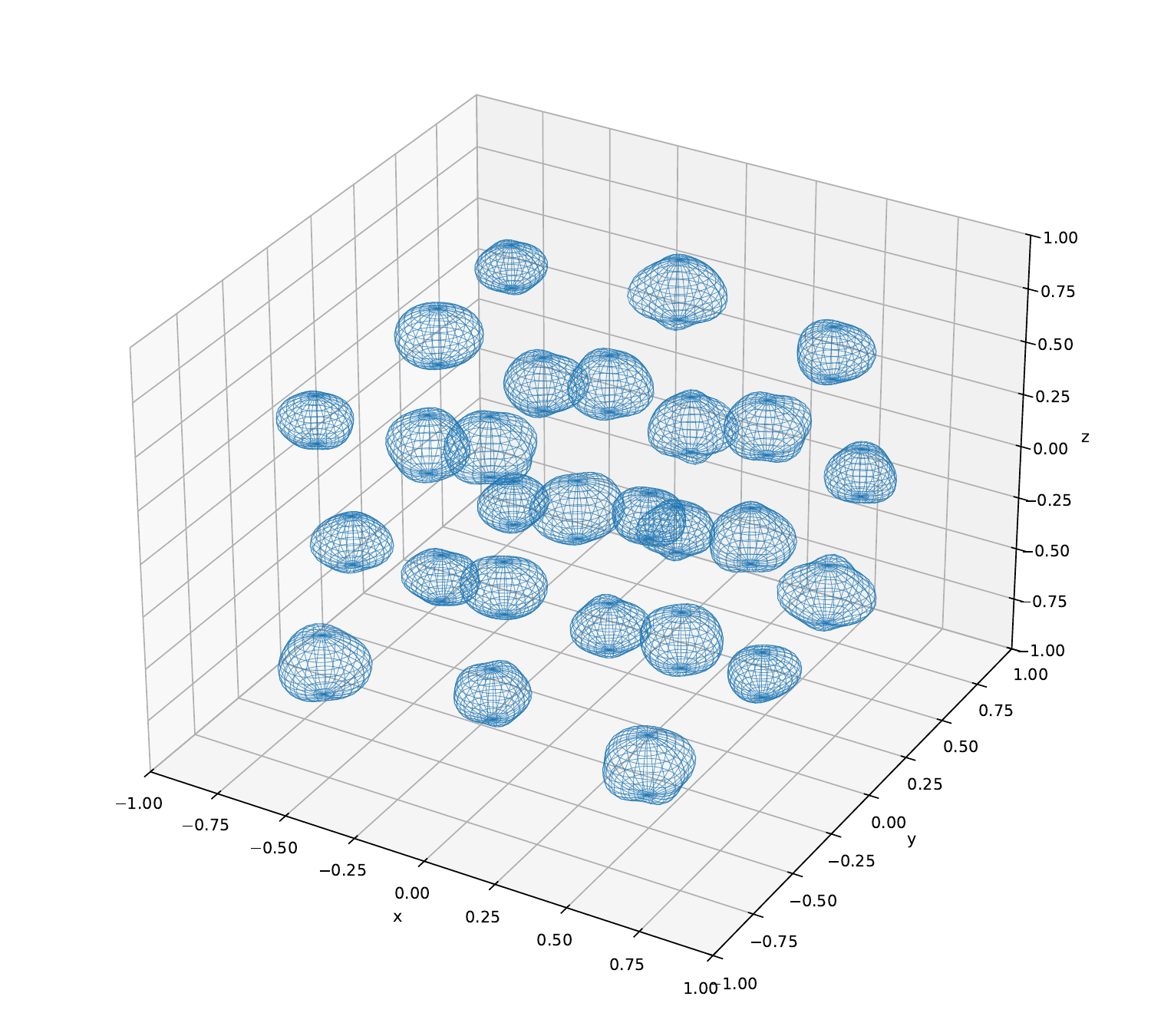}
    \caption{27 complex interfaces in 3D visualization.}
    \label{f.case4interface}
\end{figure}

Finally, we consider a 3D interface problem with 27 irregular inclusions to further evaluate the scalability of the proposed method. As shown in Fig.~\ref{f.case4interface}, the computational domain is $\Omega=[-1,1]^3$ and each interface is a surface centered at $\mathbf c_k$ and parameterized by
\begin{equation*}
\mathbf x_k(\theta,\phi)=\mathbf c_k+r_k(\theta,\phi)
\begin{pmatrix}
\sin\phi\cos\theta\\
\sin\phi\sin\theta\\
\cos\phi
\end{pmatrix},
\phi \in [0, \pi], \theta \in [0, 2\pi],
\end{equation*}
where
\begin{align*}
r_k(\theta,\phi)=r_{0,k}\bigl[&1+a_{1,k}\sin(m_{1,k}\theta+p_{1,k})\sin\phi\nonumber\\
&+a_{2,k}\cos(m_{2,k}\theta+p_{2,k})\sin^2\phi+a_{3,k}\sin(m_{3,k}\phi+p_{3,k})\bigr].
\end{align*}
The geometric parameters are varied among the interfaces to produce different irregular shapes.

The interfaces partition the domain into a background region $\Omega_0$ and 27 inclusion regions $\{\Omega_k\}_{k=1}^{27}$. The diffusion
coefficient is piecewise constant,
$$\beta(\mathbf x)=\beta_k, ~\mathbf x\in\Omega_k,~k=0,\ldots,27,$$
where $\beta_0=1$ and the coefficients $\beta_k$ are different. The background solution is chosen as
\begin{align}
u_0(x,y,z)={}&0.25\sin(\pi x)\cos(\pi y)\sin(\pi z)+0.60x-0.45y+0.35z\nonumber\\
&+0.40(x^2-y^2)+0.25xy-0.20xz+0.18yz+0.10e^{0.7x-0.5y+0.4z}.
\end{align}
Inside the $k$-th inclusion, the exact solution is constructed as
\begin{equation}
u_k(\mathbf x)=\frac{\beta_0}{\beta_k}u_0(\mathbf x)+s_k,~k=1,\ldots,27,
\end{equation}
where $s_k$ is an inclusion-dependent constant shift. 
As a result, the corresponding solution jump is
\begin{equation}
v_k(\mathbf x)=\left(\frac{\beta_0}{\beta_k}-1\right)u_0(\mathbf x)+s_k.
\end{equation}
and the source term is identical across all subdomains and is given by $f(\mathbf x)=-\beta_0\Delta u_0(\mathbf x)$.
We refer the readers to Table~\ref{T.case4_parameters} for the values of coefficients $\beta_k$ and shifts $J_k$.

\begin{table*}[!ht]
\centering
\caption{Diffusion coefficients and solution shifts for the 27 inclusions in Example~4.}
\label{T.case4_parameters}
\begin{tabular}{ccc|ccc|ccc}
\toprule[2pt]
$k$ & $\beta_k$ & $s_k$
& $k$ & $\beta_k$ & $s_k$
& $k$ & $\beta_k$ & $s_k$ \\
\midrule[1.5pt]
1 & 0.8  &  0.40 & 10 & 9.0  & -0.40 & 19 & 1.4  &  0.38 \\
2 & 8.0  & -0.25 & 11 & 2.0  &  0.50 & 20 & 7.5  & -0.33 \\
3 & 1.2  &  0.35 & 12 & 18.0 & -0.15 & 21 & 0.7  &  0.27 \\
4 & 12.0 & -0.30 & 13 & 0.6  &  0.20 & 22 & 13.0 & -0.37 \\
5 & 1.5  &  0.45 & 14 & 11.0 & -0.45 & 23 & 1.8  &  0.48 \\
6 & 10.0 & -0.20 & 15 & 2.5  &  0.32 & 24 & 20.0 & -0.18 \\
7 & 0.7  &  0.30 & 16 & 14.0 & -0.28 & 25 & 1.0  &  0.36 \\
8 & 15.0 & -0.35 & 17 & 0.9  &  0.42 & 26 & 10.5 & -0.31 \\
9 & 1.0  &  0.25 & 18 & 16.0 & -0.22 & 27 & 0.5  &  0.44 \\
\bottomrule[2pt]
\end{tabular}
\end{table*}

We sample 20000 points in $\Omega_0$ and 500 points in each inclusion region. In addition, 1000 points are sampled on each boundary face, and 512 points are generated on each interface. Table~\ref{tab:ex4} gives the accuracy and model complexity for the problem with 27 irregular interfaces. All three methods achieve comparable overall solution accuracy, with DSAN and XPINN obtaining the lowest MSE. DSAN further achieves the lowest jump and flux errors, indicating that its domain-specific representation remains effective as the number of subdomains increases. More importantly, DSAN requires only 99057 trainable parameters, compared with 1867292 for XPINN, resulting in a considerable reduction in parameter count. These results demonstrate that DSAN maintains accurate domain and interface approximations as the number of subdomains increases, while avoiding the substantial parameter growth associated with independent local networks.

\begin{table}[!ht]
\centering
\caption{Accuracy and trainable parameter comparison for the 27 irregular interfaces problem. The best result for each measurement is highlighted in \textbf{bold}.}
\label{tab:ex4}
\begin{tabular}{lcccc}
\toprule
Method & Parameters & MSE & Jump MSE & Flux MSE\\
\midrule
DSAN & 99057 & \textbf{1.49}$\times$\textbf{10}$^{-7}$ & \textbf{1.09}$\times$\textbf{10}$^{-9}$& \textbf{1.19}$\times$\textbf{10}$^{-7}$\\
CEPINN & \textbf{67001} & $1.93{\times}10^{-7}$ & $1.37{\times}10^{-9}$ & $5.75{\times}10^{-7}$ \\
XPINN   & 1867292 & \textbf{1.49}$\times$\textbf{10}$^{-7}$ & $2.60{\times}10^{-9}$ & $6.48{\times}10^{-6}$ \\
\bottomrule
\end{tabular}
\end{table}

\subsection{Scalability with Respect to the Number of Interfaces}

\begin{figure}[!ht]
    \centering
    \includegraphics[width=0.6\linewidth]{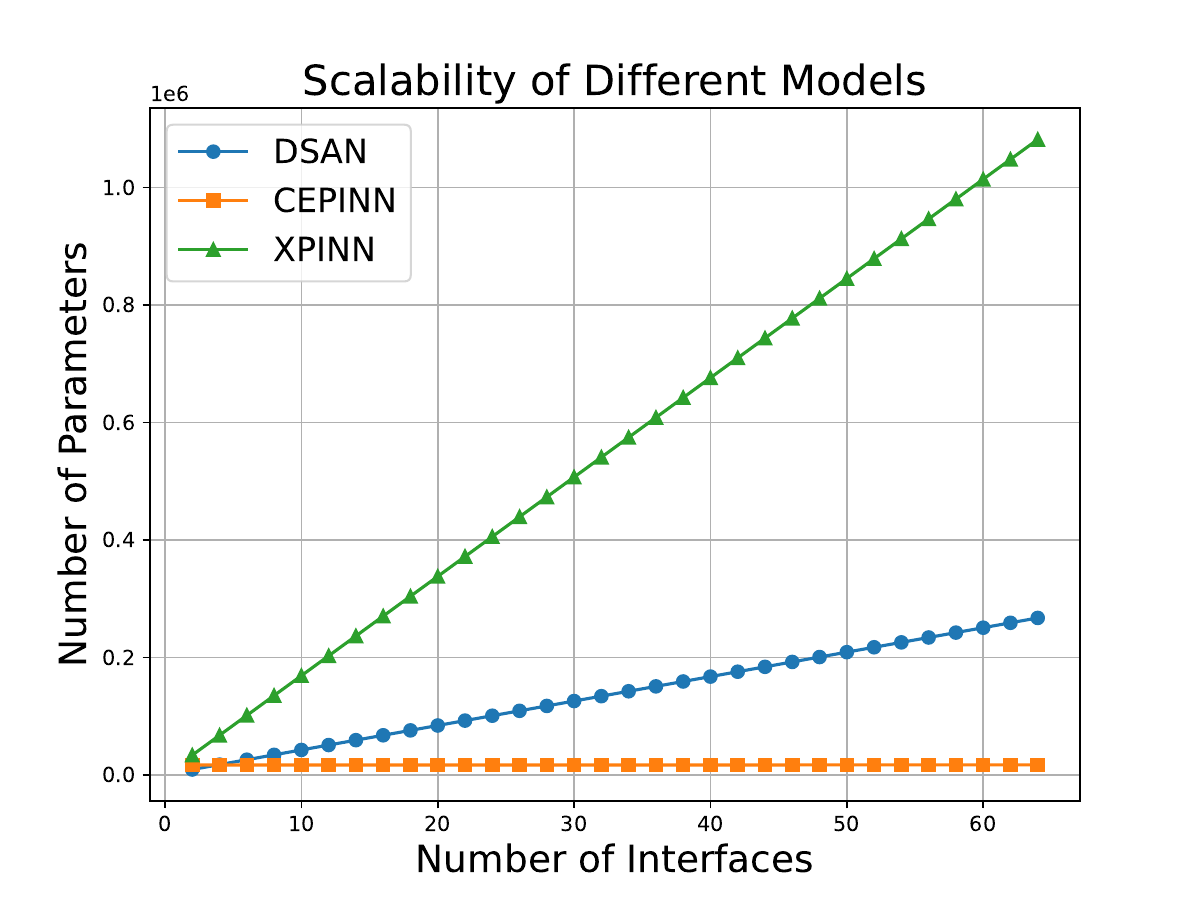}
    \caption{Number of trainable parameters of DSAN, CEPINN, and XPINN with an increasing number of interfaces.}
    \label{f.scalability}
\end{figure}

To further investigate the scalability of different network representations, we compare the number of trainable parameters of DSAN,
CEPINN, and XPINN as the number of interfaces increases from 2 to 64. The input and output dimensions are fixed at 2 and 1, respectively. The depth and width of MLP are 4 and 64, respectively. The adapter dimension for DSAN is fixed at 4, while the embedding dimension for CEPINN is set to 2. The resulting parameter counts are shown in Fig.~\ref{f.scalability}.

As the number of interfaces increases, the parameter count of XPINN grows rapidly because an independent neural network is assigned to each subdomain. In contrast, CEPINN maintains a nearly constant parameter count because subdomain information is incorporated into a single global network, with the increase in parameters arising only from the input layer. DSAN lies between these two cases, as its shared MLP remains unchanged while only lightweight adapters are added as new subdomains are introduced. Consequently, its parameter count grows much more slowly than that of XPINN. This difference becomes increasingly pronounced as the number of interfaces increases. These results demonstrate the parameter scalability of DSAN, showing that domain-specific representations can be achieved with considerably lower parameter growth than fully decomposed multi-network architectures.

\section{Conclusions}
\label{s:conclu}
In this work, we developed a Domain-Specific Adaptation Network (DSAN) for elliptic interface problems with multiple subdomains. DSAN combines a shared backbone with lightweight subdomain-specific adapters, providing local representation flexibility without replicating a complete neural network for each region. As a result, the number of trainable parameters grows much more slowly with the number of subdomains than in fully decomposed formulations.

We also established a theoretical framework connecting the physics-informed training objective to the resulting solution accuracy.
First, a residual-to-solution stability estimate shows that the continuous PDE, boundary, and interface residuals control the $L^2$ solution error under standard elliptic transmission regularity assumptions. Second, by relating the finite-sample empirical loss to its continuous population counterpart, we obtained a total error bound that separates the effects of approximation, optimization, and sampling. The approximation term was further interpreted through the shared-plus-adaptation structure of DSAN, thereby linking the network architecture directly to the error analysis. In addition, the parameter-complexity analysis quantifies the reduction in model growth relative to fully decomposed domain-specific formulations.

Numerical experiments on problems with isotropic and anisotropic coefficients, multiple interface geometries, and different interface connectivity patterns show that DSAN achieves solution accuracy comparable to or better than XPINN while requiring substantially fewer trainable parameters. This reduction becomes increasingly pronounced as the number of subdomains grows, highlighting the scalability of the proposed domain-specific representation. Overall, the results demonstrate that DSAN maintains strong approximation accuracy while effectively controlling model growth in complex multi-interface problems.

\bibliographystyle{siamplain}
\bibliography{references}

\end{document}